\documentclass[10pt,journal,final,twocolumn]{IEEEtran}
\renewcommand{\baselinestretch}{1}

\usepackage{ifthen}
\newcommand{\PaperORReport}{Report}

\usepackage{authblk}
\usepackage{amsmath}%
\usepackage{amsfonts}%
\usepackage{amssymb}%
\usepackage{psfrag,graphicx}
\usepackage{subfigure}

\usepackage[ieee]{jphmacros2e}

\makeatletter\let\NAT@parse\undefined\makeatother 

\newcommand{\R}{\mathbb{R}}

\newcommand{\N}{\mathbb{N}}

\newcommand{\Z}{\mathbf{Z}}

\newcommand{\comment}[1]{}
\newcommand{\PBcomment}[1]{\draftnote{PB: #1}}

\newcounter{rmnum}

\title{\sc\Large Mistuning-based Control Design to Improve
Closed-Loop Stability of Vehicular Platoons}

\author{Prabir~Barooah,~\IEEEmembership{Member,~IEEE,}
         Prashant G.~Mehta,~\IEEEmembership{Member,~IEEE}
        Jo\~{a}o P. Hespanha,~\IEEEmembership{Fellow, IEEE}%
\thanks{Prabir Barooah is with the Dept.~of Mechanical and Aerospace Engineering, University of Florida, Gainesville, FL 32611 (email: pbarooah@ufl.edu), Prashant G.~Mehta is with the Dept.~of Mechanical Science and Engineering,
University of Illinois, Urbana-Champaign, IL 61801(email:mehtapg@uiuc.edu), and Jo{\~a}o P. Hespanha is with the Center for Control, Dynamical Systems, and Computation, University of California, Santa Barbara, CA 93106. (email: hespanha@ece.ucsb.edu)}
\thanks{Prabir Barooah and Jo\~{a}o Hespanha's work was supported by the Institute for
Collaborative Biotechnologies through grant DAAD19-03-D-0004 from the U.S.
Army Research Office. Prashant Mehta's work was supported by the National Science Foundation by grant CMS 05-56352.}
}

\begin{document}

\maketitle

\begin{abstract}
We consider a decentralized bidirectional control of a platoon of
$N$ identical vehicles moving in a straight line. The control
objective is for each vehicle to maintain a constant velocity and
inter-vehicular separation using only the local information from
itself and its two nearest neighbors.  Each vehicle is modeled as a
double integrator. To aid the analysis, we use continuous
approximation to derive a partial differential equation~(PDE)
approximation of the discrete platoon dynamics.  The PDE model is
used to explain the progressive loss of closed-loop stability with
increasing number of vehicles, and to devise ways to combat this
loss of stability.
\medskip

If every vehicle uses the same controller, we show that the least
stable closed-loop eigenvalue approaches zero as $O(\frac{1}{N^2})$ in
the limit of a large number ($N$) of vehicles. We then show how to
ameliorate this loss of stability by small amounts of ``mistuning'',
i.e., changing the controller gains from their nominal values. We
prove that with arbitrary small amounts of mistuning, the asymptotic
behavior of the least stable closed loop eigenvalue can be improved to
${O}(\frac{1}{N})$.
%
%
All the conclusions  drawn from analysis of the PDE model are
corroborated via numerical calculations of the state-space platoon
model.
\end{abstract}

\section{Introduction}

\begin{figure}
\psfrag{xo}{$Z_0(t)$}
\psfrag{xi}{$Z_i(t)$}
\psfrag{xn}{$Z_{N+1}(t)$}
\psfrag{i}{$i$} \psfrag{1}{$1$} \psfrag{N}{$N$}
\psfrag{0}{$0$}  \psfrag{2p}{$2\pi$}
\psfrag{yi}{$y_i$} \psfrag{yi-}{$y_{i-1}$} \psfrag{yi+}{$y_{i+1}$}
\psfrag{d}{$\delta$}
\psfrag{efi}{$e_{i}^{(f)}$}
\psfrag{ebi}{$e_{i}^{(b)}$}
\begin{center}
\subfigure[A platoon with fictitious lead and follow
  vehicles.]{\includegraphics[scale = 0.4, clip = true, trim = 0in -0.5in 0in 0in]{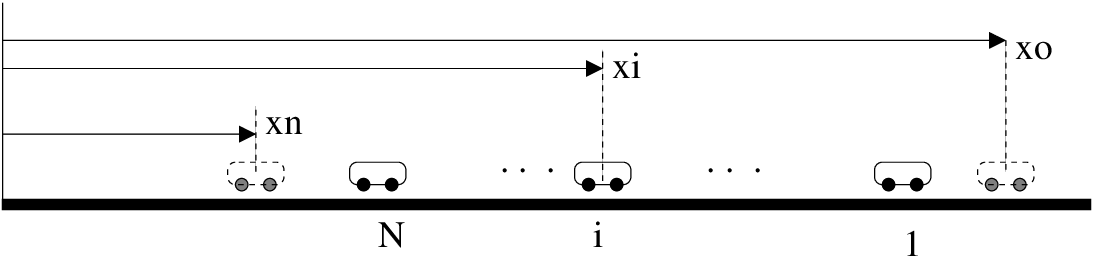}} \hspace{3 ex}
\subfigure[Same platoon in $y$ coordinates.]
{\includegraphics[scale =  0.4]{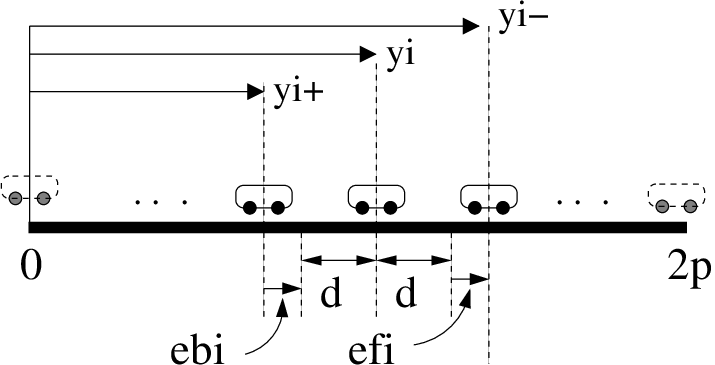}}
\end{center}
\caption{\label{fig:platoon}A platoon with $N$ vehicles moving in one dimension.}
 \end{figure}

We consider the problem of controlling a one-dimensional platoon of
$N$ identical vehicles where the individual vehicles move at a
constant pre-specified velocity $V_d$ with an inter-vehicular
spacing of $\Delta$. Figure~\ref{fig:platoon}(a) illustrates the
situation schematically. This problem is relevant to automated
highway systems (AHS) because a controlled vehicular platoon with a
constant but small inter-vehicular distance can help improve the
capacity (measured in vehicles/lane/hour, as in~\cite{DS_KH_CC_PI_VSD:94})
of a highway~\cite{JH_MT_PV_CSM:94}.
%
Due to this, the platoon control problem has been extensively
studied~\cite{RC_RH_EM_OR:58,WL_MA_TAC:66,SMM_BCK_subTAC:71,
DS_KH_CC_PI_VSD:94,MJ_BB_TAC:05,Seiler_disturb_vehicle_TAC:04}.  The
dynamic and control issues in the platoon problem are also relevant
to a general class of formation control problems including aerial
vehicles, satellites {\em
etc.}~\cite{JDW_DFC_JLS_AIAA:96,PW_FH_KL_AIAA:99}.

\medskip

Several approaches to the platoon control problem have been
considered in the literature.  These approaches fall into two broad
categories depending on the information architecture available to
the control algorithm(s): \emph{centralized} and \emph{decentralized}. We call an architecture decentralized if the control action at any individual
vehicle is computed based upon measurements obtained by on-board
sensors and possibly wireless communication with a limited
number of its neighbors. We call all other architectures centralized. Decentralized architectures investigated in the literature include the
predecessor-following~\cite{DS_KH_CC_PI_VSD:94,StStSi_TCST:00,PL_AS_TR:02} and
the bidirectional schemes~\cite{LP_TAC:74,Seiler_disturb_vehicle_TAC:04,PB_JPH_CDC:05,KC_TS:74,Zhang_IntelligentCruise_TVT:99}.
In the predecessor-following architecture, the control action at an
individual vehicle depends only on the spacing error with the
predecessor, i.e., the vehicle immediately ahead of it.  In the
bidirectional architecture, the control action depends upon relative
position measurements from both the predecessor and the follower. 
On the other hand, in a centralized architecture measurements from all the vehicles may have to be continually transmitted to a central controller or to all the
vehicles. The optimal QR designs of~\cite{WL_MA_TAC:66,MJ_BB_TAC:05}
typically lead to centralized architectures.  Predecessor and Leader
follower control schemes (see~\cite{SS_ASME:78,HT_RR_WZ_ACC:98} and
references therein), which require global information from the first
vehicle in the platoon are also examples of the centralized
architecture.
%
%
%
%
The high communication overhead in a centralized architecture makes
it less attractive for platoons with a large number of vehicles.
Additionally, with any centralized scheme, the closed loop system
becomes sensitive to communication delays that are unavoidable with
wireless communication~~\cite{XL_SM_AG_KH_ITSC:01}.

\medskip

The focus of this paper is on a decentralized bidirectional control
architecture: the control action at an individual vehicle depends
upon its own velocity and the relative position errors between
itself and its predecessor and its follower vehicles. The
decentralized bidirectional control architecture is advantageous
because, apart from its simplicity and modularity, it does not
require continual inter-vehicular communication.  Measurements
needed for the control can be obtained by on-board sensors alone.
Each vehicle is modeled as a double integrator. A double integrator
model is common in the platoon control literature since the velocity
dependent drag and other non-linear terms can usually be eliminated
by feedback linearization~\cite{DS_KH_CC_PI_VSD:94, StStSi_TCST:00}.
The control objective is to maintain a constant inter-vehicular
spacing.

\medskip

In spite of the advantages over centralized control, there are a
number of challenges in the decentralized control of a platoon,
especially when the number of vehicles, $N$, is large.  First, the
least stable closed-loop eigenvalue approaches zero as the number of
vehicles increases~\cite{PB_PM_JH_ACC:07}. Among decentralized
schemes, one particularly important special case is the so-called
{\em symmetric} bidirectional control, where all vehicles use
identical controllers that are furthermore symmetric with respect to
the predecessor and the follower position errors.  In this case, the
least stable closed loop eigenvalue approaches $0$ as
$O(\frac{1}{N^2})$ with a symmetric bidirectional control and this
behavior is independent of the choice of controller
gains~\cite{PB_PM_JH_ACC:07}.
\PBcomment{this is somewhat dissatisfactory to cite our ACC paper to
motivate the issue of loss of stability with increasing $N$, then
later prove the same result in this paper.} This progressive loss of
closed-loop damping causes the closed loop performance of the
platoon to become arbitrarily sluggish as the number of vehicles
increases.  It is interesting to note that the $O(\frac{1}{N^2})$
decay of the least stable eigenvalue occurs with the centralized LQR
control as well~\cite{MJ_BB_TAC:05}.

\medskip

The second challenge with decentralized control is that the
sensitivity of the closed loop to external disturbances increases
with increasing $N$. With predecessor following control, disturbances
acting on the vehicles cause large inter-vehicular spacing errors~\cite{RC_RH_EM_OR:58,DS_KH_CC_PI_VSD:94,SwaroopHedrick_stringstability_TAC:96}
The seminal work of~\cite{SwaroopHedrick_stringstability_TAC:96} on {\em string
instability} was partly inspired by this issue. It was shown
in~\cite{Seiler_disturb_vehicle_TAC:04} that sensitivity to
disturbances with predecessor following control is independent of
the choice of the  controller. Similar controller-independent
sensitivity to disturbances is also exhibited by the symmetric
bidirectional architecture~\cite{Seiler_disturb_vehicle_TAC:04,PB_JPH_CDC:05,RM_JB_TACsubmit:08}. In~\cite{SKY_SD_KRR_TAC:06}, it was shown that symmetric
architectures have similarly poor sensitivity even when every
vehicle uses information from more than two neighbors, as long as
the number of neighbors is no more than $O(N^{2/3})$. 


\medskip

Third, there is a lack of design methods for decentralized
architectures. For $N$ vehicles, in general, $N$ distinct
controllers need to be designed, for which few control design
methods exist. This has led to the examination of only the symmetric
control among bidirectional
architectures~\cite{Seiler_disturb_vehicle_TAC:04,PB_JPH_CDC:05,SKY_SD_KRR_TAC:06}.
Some symmetry aided simplifications are possible for analysis and
design in this case.

\medskip

In summary, while issues such as stability and sensitivity to
disturbances become critical as the platoon size increases, a lack
of analysis and control design tools in decentralized settings makes
it difficult to address these issues.

\medskip



In this paper we present a novel analysis and design method for a
decentralized bidirectional control architecture that ameliorates
the progressive loss of closed loop stability with increasing number
of vehicles. There are three contributions of this work that are
summarized below.

\medskip

First, we derive a partial differential equation (PDE) based
continuous approximation of the (spatially) discrete platoon
dynamics. Just as a PDE can be discretized using a finite difference
approximation, we carry out a reverse procedure: spatial difference
terms in the discrete model are approximated by spatial derivatives.
The resulting PDE yields the original set of ordinary differential
equations upon discretization.

\medskip

Two, we use the PDE model to derive a controller independent
conclusion on stability with symmetric bi-directional architecture.
In particular, the behavior of the least stable eigenvalue of the
discrete platoon dynamics is predicted by analyzing the eigenvalues
of the PDE.  We show that the least stable closed-loop eigenvalue
approaches zero as $O(\frac{1}{N^2})$. This prediction is confirmed
by numerical evaluation of eigenvalues for both the PDE and the
discrete platoon model.  The real part of the least stable
eigenvalue of the closed loop is taken as a measure of stability
margin.

\medskip

The third and the main contribution of the paper is a {\em
mistuning-based control design} that leads to significant
improvement in the closed loop stability margin over the symmetric
case.
The biggest advantage of using a PDE-based analysis is that the PDE
reveals, better than the state-space model does, the mechanism of
loss of stability and suggests a mistuning-based approach to
ameliorate it.  In particular, analysis of the PDE shows that
forward-backward asymmetry in the control gains is beneficial.  The
asymmetry refers to the assignment of controller gains such that a
vehicle utilizes information from the preceding and following
vehicles differently. Our main results,
Corollary~\ref{cor:mistuning-perturbation} and
Corollary~\ref{cor:mistuning-perturbation-ND}, give control gains
that achieve the best improvement in closed-loop stability by
exploiting this asymmetry.  In particular, we show that an
arbitrarily small perturbation (asymmetry) in the controller gains
from their values in the symmetric bidirectional case can result in
the least stable eigenvalue approaching $0$ only as $
{O}(\frac{1}{N})$ (as opposed to $O(\frac{1}{N^2})$ in the symmetric
bidirectional case). Numerical computations of eigenvalues of the
state-space model of the platoon is used to confirm these
predictions.  Mistuning based approaches have been used for
stability augmentation in many applications;
see~\cite{Shapiro:1998,Bendiksen:2000,RGM:2003,PM_GH_AB:07} for some
recent references.  Our paper is the first to consider such
approaches in the context of decentralized control design.

\medskip

Although the PDE model is derived under the assumption of large $N$,
in practice the predictions of the PDE model match those of the
state-space model accurately even for small values of $N$.
Similarly, the benefits of mistuning are significant even for small
values of $N$ (see Section~\ref{sec:remarks}).

\medskip

In addition to the stability margin improvements, the mistuning
design reduces the closed loop's sensitivity to external
disturbances as well. In bidirectional architectures, the $H_\infty$
norm of the transfer function from the external disturbances to the
spacing errors is used as a measure of sensitivity to disturbances~\cite{Seiler_disturb_vehicle_TAC:04}. Numerical computation of
the $H_\infty$ norm of this transfer function shows that mistuning
design also reduces sensitivity to disturbances significantly (see
Section~\ref{sec:disturbance}).

\medskip

We briefly note that there is an extensive literature on modeling
traffic dynamics using PDEs; see the seminal paper of Lighthill and
Whitham~\cite{ML_GW_PRS:55} for an early reference, the paper of
Helbing~\cite{Helbing-review:01} and references therein for a survey
of major approaches, and the papers of~\cite{DJ_CCdW_DK_CDC:05}
and~\cite{PL_RH_LA_JF_AR_TR:97} for control-oriented modeling. In
spite of apparent similarities, our approach is quite different from
the existing approaches. PDE models of traffic dynamics typically
start with continuity and momentum
equations~\cite{Helbing-review:01}. Moreover, one requires a model
of human behavior to determine an appropriate form of the external
force in the momentum equation. This difficulty frequently leads to
the introduction of terms in the PDE that are determined by fitting
data; see~\cite[Section III-D]{Helbing-review:01} for a thorough
discussion of such approximations used in various continuum traffic
models. In contrast, we approximate the closed loop dynamic
equations by continuous functions of space (and time) that are 
inspired by finite-difference discretization of PDEs. Ad-hoc
approximations of human behavior is not needed. Moreover, the
original dynamics can be recovered by discretizing the derived PDE,
which provides further evidence of consistency between the
(spatially) discrete and continuous models.

%
We also note that macroscopic models of traffic flow models
have been used for designing control laws for a complete automated
highway system (AHS) with lane changing, merging, etc. in addition to a platoon in one lane (see~\cite{PL_RH_LA_JF_AR_TR:97,LA_RH_PL_CEP:99} and references
therein). The PDE model derived in the paper is not applicable to a
complete AHS, but only to a single platoon.

\medskip

The rest of the paper is organized as follows:
section~\ref{sec:problem} states the platoon problem in formal terms
by describing a state-space model of the closed loop platoon
dynamics; section~\ref{sec:continuous-model} then describes the
derivation of the PDE model from the state space model. In
section~\ref{sec:instability-analysis} the PDE is analyzed to
explain the loss of stability with $N$, and
section~\ref{sec:mistuning} describes how to ameliorate such loss of
stability by mistuning. Section~\ref{sec:sim} reports simulation
results that show the benefit of mistuning in time-domain. In
Section~\ref{sec:remarks}, we comment on various aspects of the
proposed mistuning design.

\section{Closed loop dynamics with bidirectional control}\label{sec:problem}
Consider a platoon of $N$ identical vehicles moving in a straight line
as shown schematically in Figure~\ref{fig:platoon}(a).  Let $Z_i(t)$
and $V_i(t) \eqdef  \dot{Z}_i(t)$ denote the position and the velocity, respectively, of
the $i^{\text{th}}$ vehicle for $i=1,2,\dots,N$.  Each vehicle is modeled as a double integrator:
\begin{align}\label{eq:vehicle-model}
\ddot{Z}_i = U_i,
\end{align}
where $U_i$ is the control (engine torque) applied on the
$i^{\text{th}}$ vehicle. Such a model arises after the velocity dependent drag and
other non-linear terms have been eliminated by using feedback
linearization~\cite{DS_KH_CC_PI_VSD:94,StStSi_TCST:00}.

The control objective is to maintain a constant inter-vehicular
distance $\Delta$ and a constant velocity $V_d$ for every vehicle.
Every vehicle is assumed to know the desired spacing $\Delta$ and
the desired velocity $V_d$. The control architecture is required to
be decentralized, so that every vehicle uses locally available
measurements. We assume that the error between the position (as well
as velocity) of a vehicle and its desired value is small, so that
analysis of the platoon dynamics with a linear vehicle model and
 a linear control law is justified.

\medskip

\begin{table}
\begin{center}
\begin{tabular}{c|c|c|c}
Scenario  & Length $L$ & Leader & Follower \\ \hline &&&\\
I  & $(N+1)\Delta$  & $\tilde{v}_{0}=0$ & $\tilde{v}_{N+1}=0$ \\ \hline &&& \\
II & $N\Delta$  &  $\tilde{v}_{0}=0$ & --   \\
\hline
\end{tabular}
\end{center}
\caption{\label{tab:scenario}The two scenarios.}
\end{table}

In this paper, we assume a bi-directional control architecture for
individual vehicles in the platoon (except the first and the last
vehicles). For the first and the last vehicles, we consider two
types of control architectures (termed as scenarios I and II) as
tabulated in Table~\ref{tab:scenario}. In scenario~I, we introduce
(after~\cite{MJ_BB_TAC:05,SMM_BCK_subTAC:71}) a fictitious lead
vehicle and a fictitious follow vehicle, indexed as $0$ and $N+1$
respectively.  Their behavior is specified by imposing a constant
velocity trajectories as $Z_0(t) = V_d \, t$ and $Z_{N+1} = V_d \, t -
(N+1)\Delta$.
In scenario~II, only a fictitious lead vehicle with index $i=0$ with
$Z_0(t) = V_d t$ is introduced.  For the last vehicle in the platoon
in scenario II, there is no follower vehicle and it uses information
only from its predecessor to maintain a constant gap.

\medskip
Consistent with the decentralized bidirectional linear control architecture,
the control $U_i$ for the $i^\text{th}$ vehicle is assumed to depend
only on 1) its velocity error $V_i-V_d$, and 2) the relative position
errors between itself and its immediate neighbors. That is,
\begin{multline}
U_i  = k^{(f)}_{i}(Z_{i-1}-Z_i - \Delta) - k^{(b)}_{i}(Z_i -
Z_{i+1}-\Delta)  \\
- b_i(V_i-V_d).\label{eq:control-law-inline-1}
\end{multline}
where $k^{(\cdot)}_{i},b_i$ are positive constants. The first two
terms are used to compensate for any deviation away from nominal
position with the predecessor (front) and the follower (back) vehicles
respectively.  The superscripts $^{(f)}$ and $^{(b)}$ correspond to
\emph{front}  and  \emph{back}, respectively. The third term is used
to obtain a zero steady-state error in velocity. In principle,
relative velocity errors between neighboring vehicles can also be
incorporated into the control, but we do not examine this situation
here. Since $V_d$ and $\Delta$ are known to every vehicle, the
relative errors used in the control law, including the velocity
error, can be obtained in practice by on-board devices such as
radars, GPS, and speed sensors.

\medskip

The control law~\eqref{eq:control-law-inline-1} represents state
feedback with local (nearest neighbor) information. Analysis of
this controller structure is relevant even if there are additional
dynamic elements in the controller. There are several reasons for
this.  First, a dynamic controller cannot have a zero at the origin.
It will result in a pole-zero cancellation causing the steady-state
errors to grow without bound as $N$ increases~\cite{PB_JPH_CDC:05}.
Second, a dynamic controller cannot have an integrator either. For
if it does, the closed-loop platoon dynamics become unstable for a
sufficiently large values of $N$~\cite{PB_JPH_CDC:05}.  As a result,
any allowable dynamic compensator must essentially act as a static
gain at low frequencies.  The results of~\cite{PB_JPH_CDC:05}
indicate that the principal challenge in controlling large platoons
arises due to the presence of a double integrator with its unbounded
gain at low frequencies. Hence, the limitation and its amelioration
discussed here with the local state feedback structure
of~(\ref{eq:control-law-inline}) is also relevant to the case where
additional dynamic elements appear in the control.

\medskip

To facilitate analysis, we consider a coordinate change
\begin{equation}
 y_i  = 2\pi(\frac{Z_i(t) - V_d t+ L}{L}), \qquad
 v_i  = 2\pi\frac{V_i - V_d }{L},
\label{eq:coordinate-change}
 \end{equation}
where $L$ denotes the \emph{desired platoon length}, which equals $(N+1)\Delta$ in scenario~I and $N\Delta$ in scenario~II. Figure~\ref{fig:platoon}(b)
depicts the schematic of the platoon in the new
coordinates. The scaling ensures that $y_0(t) \equiv 2\pi$,
$y_i(t)\in [0,2\pi]$, and $y_{N+1}(t) \equiv 0$ ($y_N(t)=0$) in
scenario~I (II).  Here, we have implicitly assumed that deviations of
the vehicle positions and velocities from their desired values are
small.

\medskip

In the scaled coordinate, the dynamics of the $i^\text{th}$ vehicle
are described by
\begin{align}\label{eq:vehicle-dynamics-y}
\ddot{y}_i = u_i,
\end{align}
where $u_i \eqdef 2\pi U_i/L$.  The
desired spacing and velocities are
\begin{align}\label{eq:delta-defn}
\delta := \frac{\Delta}{L/2\pi}, \qquad v_d:= \frac{V_d -
V_d}{L/2\pi}=0,
\end{align}
and the desired position of the $i^{\text{th}}$ vehicle is
\begin{align}\label{eq:ydesired}
y_i^{d}(t) \equiv 2\pi - i \delta.
\end{align}
The position and velocity errors for the $i^{\text{th}}$ vehicle are
given by:
\begin{align}\label{eq:y-v-tilde-defn}
\begin{split}
\tilde{y}_i(t) & = y_i(t) - y_i^d(t),  \tilde{v}_i  = v_i -
v_d = v_i, \text{and } \\
\dot{\tilde{y}}_i & = \tilde{v}_i. 
\end{split}
\end{align}
We note that $\tilde{v}_0 = \tilde{v}_{N+1} = 0$ for the fictitious
lead and follow vehicles. In the scaled coordinates, the decentralized bidirectional control law~\eqref{eq:control-law-inline-1} is equivalent to the following
\begin{align}
u_i &= k^{(f)}_{i}\,( y_{i-1} - y_i - \delta) - k^{(b)}_{i}\,(y_{i} - y_{i+1} -\delta)   - b_i\,\tilde{v}_i \label{eq:control-law-inline}\\
& = k^{(f)}_{i}\,(\tilde{ y}_{i-1} - \tilde{y}_i)  - k^{(b)}_{i}\,(\tilde{y}_{i} - \tilde{y}_{i+1}) - b_i \tilde{v}_i. \nonumber
\end{align}
It follows from~\eqref{eq:vehicle-dynamics-y}
and~\eqref{eq:control-law-inline} that the closed loop dynamics of
the $i^\text{th}$ vehicle in the $\tilde{y}$-coordinate is
\begin{align}\label{eq:closed-loop-y-i-final}
\ddot{\tilde{y}}_i + b_i \dot{\tilde{y}}_i = k_i^{(f)}(\tilde{y}_{i-1}- \tilde{y}_i) - k_i^{(b)}(\tilde{y}_i - \tilde{y}_{i+1}).
\end{align}
To describe the  closed-loop dynamics of the whole platoon, we define
\begin{eqnarray}
\mbf{\tilde
y}\eqdef [\tilde{y}_1,\tilde{y}_2,\dots,\tilde{y}_N]^T, \qquad \mbf{\tilde {v}}\eqdef [\tilde{v}_1,\dots,\tilde{v}_N]^T.\nonumber
\end{eqnarray}
For scenario~I with fictitious lead and follow vehicles, the control law~\eqref{eq:control-law-inline} yields the following closed loop dynamics.
\begin{align}\label{eq:closedloop-platoon-fb} \matt{\mbf{\dot
{\tilde y}} \\ \mbf{\dot{\tilde v}}} & = \underbrace{\matt{ 0 & I
\\-K^{(f)}_{\mathrm I} M^T - K^{(b)}_\mathrm{I} M & -B}}_{A_{L-F}}\matt{\mbf{\tilde y} \\ \mbf{\tilde v}}
\end{align}
where $K^{(f)}_\mathrm{I} =\diag(k_{1}^{(f)},k_{2}^{(f)},\dots,k_{N}^{(f)})$,
$K^{(b)}_\mathrm{I} = \diag(k_{1}^{(b)},k_{2}^{(b)},\dots,k_{N}^{(b)})$, $B =
\diag(b_{1},b_{2},\dots,b_{N})$,  and
\begin{align*}
M & = \smatt{1 & -1 & 0 & \dots & &\\ 0 & 1 & -1 & & & \\ \vdots & &
& \ddots & & 0\\ & & & & 1 & -1\\& & & \dots & 0 & 1}.
\end{align*}
For scenario~II with a fictitious lead vehicle and no follow vehicle, the closed loop
dynamics are
\begin{align}\label{eq:closedloop-platoon-f}
\matt{\mbf{\dot {\tilde y}} \\ \mbf{\dot{\tilde v}}} & =
\underbrace{ \matt{ 0 & I
\\-K^{(f)}_\mathrm{II} M^T - K^{(b)}_\mathrm{II} M_o  & -B}}_{A_L}\matt{\mbf{\tilde y} \\ \mbf{\tilde v}},
\end{align}
where $K^{(f)}_\mathrm{II} = K^{(f)}_\mathrm{I}$, $K^{(b)}_\mathrm{II} =
\diag(k_{1}^{(b)},k_{2}^{(b)},\dots,k_{N-1}^{(b)},0)$, and
\begin{align*}
M_o & = \smatt{1 & -1 & 0 & \dots & &\\ 0 & 1 & -1 & & & \\ \vdots &
& & \ddots & & 0\\ & & & & 1 & -1\\& & & \dots & 0 & 0}.
\end{align*}

\medskip

Our goal is to understand the behavior of the closed loop stability
margin with increasing $N$ and to devise ways to improve it by
appropriately choosing the controller gains. While in principle this
can be done by analyzing the eigenvalues of the matrix $A_{L-F}$
(scenario~I) and of $A_L$ (scenario~II), we take an alternate route.
For large values of $N$, we approximate the dynamics of the discrete
platoon by a partial differential equation (PDE) which is used for
analysis and control design.

\section{PDE model of platoon closed loop dynamics}\label{sec:continuous-model}
In this section, we develop a continuous PDE approximation of the
(spatially) discrete platoon dynamics.  The PDE is derived with
respect to a scaled spatial coordinate $x\in [0,2\pi]$. We recall
that in Section~\ref{sec:problem}, the scaled location of the
$i^{\text{th}}$ vehicle (denoted as $y_i$) was defined with
respect to such a coordinate system. In effect, the two symbols $x$
and $y$ correspond to the same coordinate representation but are
used here to distinguish the continuous and discrete formulations.
As in the discrete case, the platoon
always occupies a length of $2\pi$ irrespective of $N$.

\subsection{PDE derivation}
The starting point is a continuous approximation:
\begin{align*}
v(x,t) & \eqdef v_i(t) \quad \text{ at } x = y_i\\
\Rightarrow v(x,t) & = \tilde{v}_i(t). \quad \text{ (from~\eqref{eq:y-v-tilde-defn})}
\end{align*}
Similarly, $b(x), k^{(f)}(x), k^{(b)}(x)$ are used to denote continuous
approximations of discrete gains $b_i,k_i^{(f)},k_i^{(b)}$ respectively.
We will construct a PDE approximation of discrete dynamics in terms
of these continuous approximations.
%
To do so, it is convenient to first
differentiate~\eqref{eq:closed-loop-y-i-final} with respect to time,
\begin{align}\label{eq:dynamics-i-in-velocity}
\ddot{\tilde{v}}_i + b_i \dot{\tilde{v}}_i = k_i^{(f)}(\tilde{v}_{i-1}-\tilde{v}_i) - k_i^{(b)}(\tilde{v}_i - \tilde{v}_{i+1}).
\end{align}
We recast this equation
\begin{multline*}
\ddot{\tilde{v}}_i  + b_i \dot{\tilde{v}}_i  = - k_i^{(+)}\tilde{v}_i + \frac{1}{2}(k_i^{(+)} + k_i^{(-)})\tilde{v}_{i-1}
 - \frac{1}{2}( k_i^{(+)} - k_i^{(-)})\tilde{v}_{i+1},
\end{multline*}
where
\begin{align}\label{eq:kplusminus-def}
k_i^{(+)} & \eqdef k_i^{(f)} + k_i^{(b)}, & k_i^{(-)} & \eqdef k_i^{(f)} - k_i^{(b)}.
\end{align}
It follows that
\begin{align*}
\ddot{\tilde{v}}_i  + b_i \dot{\tilde{v}}_i & = \frac{1}{2}k_i^{(-)}(\tilde{v}_{i-1} - \tilde{v}_{i+1}) + \frac{1}{2}k_i^{(+)}(\tilde{v}_{i-1} - 2\tilde{v}_i + \tilde{v}_{i+1}) \\
& =  \frac{1}{\rho_0}k_i^{(-)}\frac{\tilde{v}_{i-1} - \tilde{v}_{i+1}}{2\delta} + \frac{1}{2\rho_0^2}k_i^{(+)}\frac{\tilde{v}_{i-1} - 2\tilde{v}_i + \tilde{v}_{i+1}}{\delta_0^2}
\end{align*}
where
\begin{align}
\rho_0 \eqdef \frac{1}{\delta} = \frac{N}{2\pi}.
\end{align}
$\rho_0$ has the physical interpretation of the \emph{mean density}
(vehicles per unit length). Now, we make a finite-difference
approximation of derivatives
%
\begin{align*}
\frac{\tilde{v}_{i-1} - \tilde{v}_{i+1}}{2\delta} = \left[ \frac{\partial}{\partial x} v(x,t)\right]_{x = y_i} \\
\frac{\tilde{v}_{i-1} - 2\tilde{v}_i + \tilde{v}_{i+1}}{\delta_0^2} = \left[ \frac{\partial^2}{\partial x^2} v(x,t)\right]_{x = y_i},
\end{align*}
where we recall that $v(x,t)$ is a continuous approximation of the
vehicle velocities ($\tilde{v}_i(t) = v(y_i,t)$ etc).  Denoting $k^{(+)}(x)$
and $k^{(-)}(x)$ as continuous approximations of $k_i^{(+)}$ and
$k_i^{(-)}$ respectively, the discrete model is written as:
\begin{multline*}
\left[\frac{\partial^2}{\partial t^2}v(x,t)\right]_{x = y_i}+ \left[ b(x) \frac{\partial}{\partial t}v(x,t) \right]_{x = y_i} = \\
\frac{1}{\rho_0} \left[ k^{(-)}(x) \frac{\partial  }{\partial x}v(x,t)\right]_{x = y_i} + \frac{1}{2\rho_0^2} \left[ k^{(+)}(x)\frac{\partial^2 }{\partial x^2}v(x,t) \right]_{x = y_i}
\end{multline*}
Hence, we arrive at the partial differential equation (PDE) as a
model of the discrete platoon dynamics:
\begin{multline}\label{eq:PDE-v-1}
\left(\frac{\partial^2}{\partial t^2} + b(x)\frac{\partial}{\partial t} \right)v(x,t) =\\
 \left( \frac{1}{\rho_0} k^{(-)}(x) \frac{\partial }{\partial x} + \frac{1}{2\rho^2_0} k^{(+)}(x) \frac{\partial^2 }{\partial x^2}\right)v(x,t)
\end{multline}
In the remainder of this paper, we assume that $k^{(+})(x)>0$.
Using~\eqref{eq:kplusminus-def}, the continuous counterparts of the
front and the back gains are given by
\begin{align}\label{eq:kfrontback-x-def}
\begin{split}
k^{(f)}(x) & = \frac{1}{2}\left(k^{(+)}(x) + k^{(-)}(x)\right), \\
k^{(b)}(x) & = \frac{1}{2}\left(k^{(+)}(x) - k^{(-)}(x)\right),
\end{split}
\end{align}
so that the gain values $k_i^{(\cdot)}$ can be obtained as
$k_i^{(f)} = k^{(f)}(y_i)$ and  $k_i^{(b)} = k^{(b)}(y_i)$.  It can
be readily verified that one recovers the system of ordinary
differential equations (\eqref{eq:dynamics-i-in-velocity} for
$i=1,\dots,N$) by discretizing the PDE~\eqref{eq:PDE-v-1} using a
finite difference scheme on the interval $[0, 2\pi]$ with a
discretization $\delta$ between discrete points.

\medskip

The boundary conditions for the PDE~\eqref{eq:PDE-v-1} depend upon the dynamics of the
first and the last vehicles in the platoon.  For scenario~I with a
constant velocity fictitious lead and follow vehicles, the appropriate
boundary conditions are of the Dirichlet type on both ends:
\begin{equation}\label{eq:Dirichlet-bc}
 v (0,t)= v (2\pi,t)=0, \quad \forall t\in[0,\infty).
\end{equation}
For scenario~II with the only a fictitious lead vehicle, the
appropriate
 boundary conditions are of Neumann-Dirichlet type:
\begin{align}\label{eq:Neumann-Dirichlet-bc}
\frac{\partial v}{\partial x} (0,t) & = v(2\pi,t) = 0. \quad \forall t\in[0,\infty)
\end{align}
We refer the reader to Appendix~\ref{sec:pde-properties} for a
discussion on well-posedness of the solutions to~\eqref{eq:PDE-v-1}.
It is shown in  Appendix~\ref{sec:pde-properties} that a solution
exists in a weak sense when $k^{(+)},k^{(-)},\frac{d
k^{(+)}}{dx} \in L^\infty([0,2\pi])$.

\smallskip

Equation~\eqref{eq:PDE-v-1} describes spatio-temporal evolution of small velocity perturbations in a platoon.  It is worthwhile to note that the PDE
model is a hyperbolic equation.  Without the two first order terms
(i.e., for $b(x)=k^{(-)}(x)=0$), the PDE is a standard wave equation
with spatially inhomogeneous values of wave speed.  The term
$\frac{1}{\rho_0}k^{(-)}(x)\frac{\partial v}{\partial x}$ is an
advection term, and $b(x)\frac{\partial v}{\partial t}$ is a damping
term.  The hyperbolic nature of the PDE model means that a
perturbation originating, say, in the middle of a long platoon will
propagate both upstream and downstream with finite speed.  The two
first order terms serve to modify aspects of this propagation.  The
damping term causes a perturbation to damp out in time.  The
advection term serves to create possible asymmetries in upstream
versus downstream propagation.

\subsection{Eigenvalue comparison}
For preliminary comparison of the PDE obtained above with the
state-space model of the closed loop platoon dynamics, we consider
the simplest case where the position control gains are constant for
every vehicle, i.e., $k^{(f)}(x) = k^{(b)}(x) = k_{0}$ and
$b(x)=b_0$.  In such a case $k^{(-)}(x) \equiv 0$, $k^{(+)}(x)
\equiv 2k_0$ and the PDE~\eqref{eq:PDE-v-1} simplifies to
\begin{align}\label{eq:PDE-rho-2}
 \left( \frac{\partial^2}{\partial t^2} + b_0\frac{\partial}{\partial
 t} - \frac{k_0}{\rho_0^2}\frac{\partial^2}{\partial x^2}
 \right) v = 0,
\end{align}
which is a damped wave equation with a wave speed of
$\frac{\sqrt{k_0}}{\rho_0}$.  The wave equation is consistent with
the physical intuition that a symmetric bidirectional control
architecture causes a disturbance to propagate equally in both
directions.

\begin{figure}
\begin{center}
\subfigure[Scenario I ( Dirichlet-Dirichlet)]{\includegraphics[scale=0.3]{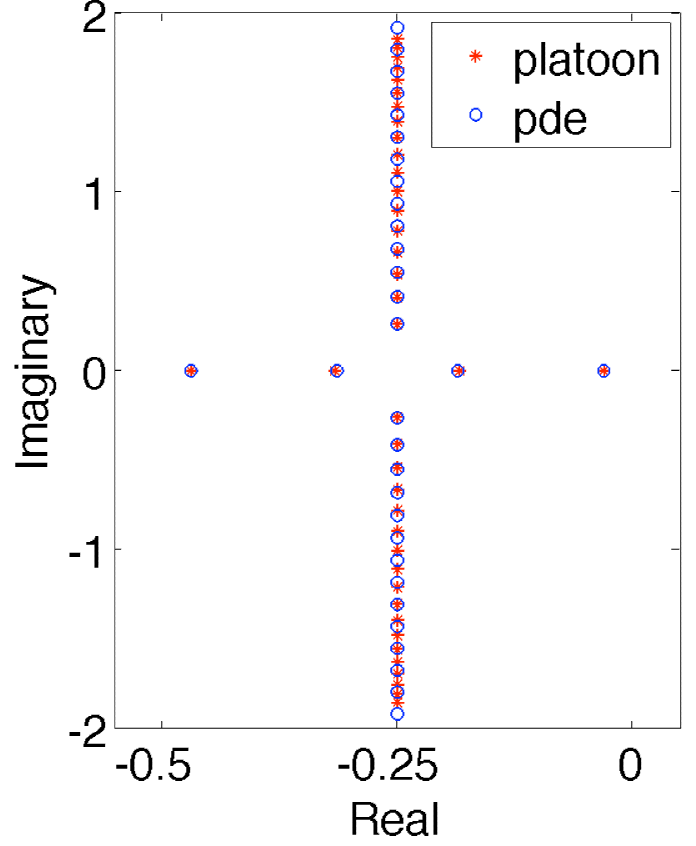}}
\hspace{6 ex} \subfigure[Scenario II ( Neumann-Dirichlet)]{\includegraphics[scale=0.3]{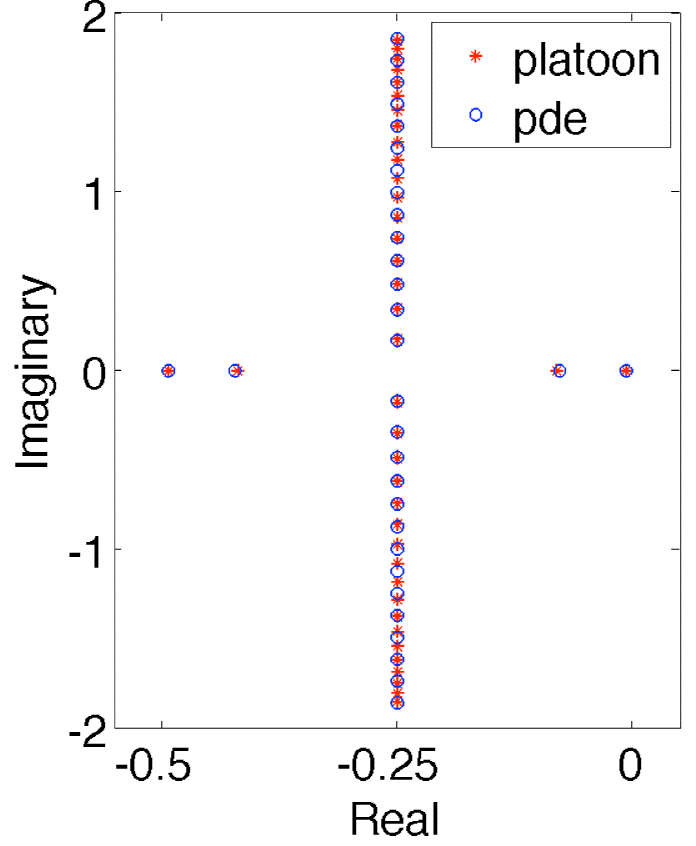}}
\end{center}
\caption{\label{fig:PDEvalidation-spectra1} Comparison of closed loop
eigenvalues of the platoon dynamics and the eigenvalues of the
corresponding PDE~\eqref{eq:PDE-rho-2} for the two different
scenarios: (a) platoon with fictitious lead and follow vehicles,
and correspondingly the PDE~\eqref{eq:PDE-rho-2} with Dirichlet
boundary conditions, (b) platoon with fictitious lead vehicle, and
correspondingly the PDE~\eqref{eq:PDE-rho-2} with
Neumann-Dirichlet boundary conditions. For ease of comparison, only
a few of the eigenvalues are shown. Both plots are for $N=25$
vehicles; the controller parameters are $k^{(f)}_i=k^{(b)}_i = 1$
and $b_i=0.5$ for $i=1,2,\dots,N$, and for the PDE $k^{(f)}(x)\equiv k^{(b)}(x) \equiv 1$ and $b(x)\equiv 0.5$.}
\end{figure}

Figure~\ref{fig:PDEvalidation-spectra1} compares the closed loop
eigenvalues of a discrete platoon with $N=25$ vehicles and the
PDE~\eqref{eq:PDE-rho-2}.  The eigenvalues of the platoon are
obtained by numerically evaluating the eigenvalues of the matrices
$A_{L-F}$ and $A_L$ (defined in ~\eqref{eq:closedloop-platoon-fb}
and \eqref{eq:closedloop-platoon-f}).  The eigenvalues of the PDE
are computed numerically after using a Galerkin method with Fourier
basis~\cite{Canuto}.  The figure shows that the two sets of
eigenvalues are in  excellent match. In particular, the least stable
eigenvalues are well-captured by the PDE. Additional comparison
appears in the following sections, where we present the results for
analysis and control design.

\section{Analysis of the symmetric bidirectional case}\label{sec:instability-analysis}

This section is concerned with asymptotic formulas for stability
margin (least stable eigenvalue) for the symmetric bidirectional
architecture with symmetric and constant control gains: $k^{(f)}(x)
= k^{(b)}(x) \equiv k_{0}$ and $b(x) \equiv b_0$.  The analysis is
carried out with the aid of the associated PDE model:
\begin{equation}\label{eq:PDE-nominal}
\left(\frac{\partial^2}{\partial t^2} + b_0\frac{\partial}{\partial t}- a_0^2\frac{\partial^2}{\partial x^2}\right)\tilde v=0,
\end{equation}
where $x\in [0,2\pi]$ and
\begin{align}\label{eq:a0-def}
a_0^2 : = \frac{k_0}{\rho_0^2}
\end{align}
is the wave speed.  The closed-loop eigenvalues of the PDE
require consideration of the eigenvalue problem
\begin{align}\label{eq:PDE-eigenproblem}
\frac{d^2 \eta}{dx^2} = \lambda \eta(x),
\end{align}
where $\eta$ is an eigenfunction that satisfies appropriate boundary
conditions:~\eqref{eq:Dirichlet-bc} for scenario~I
and~\eqref{eq:Neumann-Dirichlet-bc} for scenario~II.  The
eigensolutions to the eigenvalue problem~\eqref{eq:linPDE_lap} for
the two scenarios are given in Table~\ref{tab:eigen-laplacian}.  The
eigenfunctions in either scenario provide a basis of
$L^2([0,2\pi])$.

\begin{table}
\begin{center}
\begin{tabular}{c|c|c|c}
boundary condition & eigenvalue $\lambda_l$ & \begin{minipage}{.12\columnwidth} eigenfunction \\
$\psi_l(x)$ \end{minipage} & $l$ \\ \hline &&& \\
\begin{minipage}{.14\textwidth}$\eta(0)=\eta(2\pi) = 0$ \\ (Dirichlet - Dirichlet) \\\end{minipage} & $-\frac{l^2}{4}$  & $\sin(\frac{l x}{2})$ & $l=1,2,\dots$ \\ \hline &&& \\
\begin{minipage}{.14\textwidth} $\frac{\partial \eta}{\partial x}(0) = \eta(2\pi) = 0$ \\   (Neumann - Dirichlet) \\\end{minipage} & $-\frac{(2l-1)^2}{16}$  & $\cos(\frac{(2l - 1)x}{4})$ & $l=1,2,\dots$ \\ \hline
\end{tabular}
\vspace{1 ex}
\end{center}
\caption{\label{tab:eigen-laplacian}The eigen-solutions for the
Laplacian operator with two different boundary conditions.}
\end{table}

\medskip
%

After taking a Laplace transform, the eigenvalues of the PDE
model~\eqref{eq:PDE-nominal} are obtained as roots of the
characteristic equation
\begin{equation}\label{eq:linPDE_lap}
s^2 + b_0 s - a_0^2 \lambda = 0,
\end{equation}
where $\lambda$ satisfies~\eqref{eq:PDE-eigenproblem}.  Using
Table~\ref{tab:eigen-laplacian}, these roots are easily evaluated.
For instance, the $l^{\text{th}}$ eigenvalue of the PDE~\eqref{eq:PDE-nominal} with Dirichlet boundary conditions is given by
\begin{equation}
s_l^{\pm}=\frac{-b_0\pm\sqrt{b_0^2- a_0^2 l^2}}{2},\label{eq:sksol}
\end{equation}
where $l=1,2,\hdots$.  The real part of the eigenvalue depends upon
the discriminant $D(l,N) \eqdef (b_0^2-a_0^2 l^2)$, where the wave speed
$a_0$ depends both on control gain $k_0$ and number of vehicles $N$
(see~(\ref{eq:a0-def})).  For a fixed control gain, there are two
cases to consider:
\begin{enumerate}
\item If $D(l,N)<0$, the roots $s_l^{\pm}$ are complex with the
real part given by $-\frac{b_0}{2}$,
\item If $D(l,N)>0$, the roots $s_l^{\pm}$ are real with
$s_l^{+} + s_l^{-} = - b_0$.
\end{enumerate}
In the former case, the damping is determined by the velocity
feedback term $b_0\frac{\partial }{\partial t}$, while in the latter
case one eigenvalue ($s_l^{-}$) gains damping at the expense of the
other ($s_l^{+}$) which looses damping.  When $s_l^{\pm}$ are real,
the eigenvalue $s_l^+$ is closer to the origin than $s_l^{-}$; so we
call $s_l^+$ the $l^\text{th}$ \emph{less-stable} eigenvalue.
The following lemma gives the asymptotic formula for this eigenvalue
in the limit of large $N$.


\medskip

\begin{lemma}\label{lem:s-ell-symmetric}
Consider the eigenvalue problem for the linear
PDE~\eqref{eq:PDE-nominal} with boundary
conditions~(\ref{eq:Dirichlet-bc})
and~(\ref{eq:Neumann-Dirichlet-bc}), corresponding to scenarios~I
and~II respectively.  The $l^{\text{th}}$ less-stable eigenvalue
$s^+_l$ approaches $0$ as $O(1/N^2)$ in the limit as
$N\rightarrow\infty$. The asymptotic formulas appear in
Table~\ref{tab:eigenvalue-asymptotes-nominal}.\frqed
\end{lemma}

\begin{table}
\begin{center}
\begin{tabular}{l|c|c}
boundary condition & $s_l^+$ for $l<< l_c$ & $l_c$\\ \hline && \\ Dirichlet-Dirichlet & $ - \frac{\pi^2k_0}{b_0}\frac{l^2}{N^2}+ O(\frac{1}{N^4})$ & $\frac{b_0 N}{2\pi\sqrt{k_0}}$ \\  && \\
\hline & & \\
Neumann-Dirichlet & $ -\frac{\pi^2 k_0}{4
b_0}\frac{l^2}{N^2}+O(\frac{1}{N^4})$  & $\frac{b_0
N}{2\pi\sqrt{k_0}}$
\end{tabular}\end{center}
\caption{\label{tab:eigenvalue-asymptotes-nominal}the trend of the
less stable eigenvalue $s_l^{+}$ for the PDE~\eqref{eq:PDE-nominal}}
\end{table}

\begin{proof-lemma}{\ref{lem:s-ell-symmetric}}
We first consider scenario~I with Dirichlet boundary
conditions~(\ref{eq:Dirichlet-bc}). Using~(\ref{eq:sksol})
and~(\ref{eq:a0-def}),
\begin{align*}
2s_l^{\pm} & =-b_0 \pm b_0 \left(1-\frac{a_0^2 l^2}{b_0^2}\right)^{1/2} \\
 & = -b_0\pm b_0 \left(1 - \frac{2\pi^2
k_0}{b_0^2}\frac{l^2}{N^2}\right) + O(\frac{1}{N^4})
\end{align*}
for $a_0^2 l^2/b_0^2 << 1$.  The asymptotic formula holds for wave
numbers
\begin{align}\label{eq:l-critical}
l \ll \frac{b_0}{a_0} = \frac{b_0 N}{2\pi \sqrt{k_0}} \defeq l_c,
\end{align}
and in particular for each $l$ as $N\rightarrow\infty$.  The proof
for the scenario~II with Neumann-Dirichlet boundary
conditions~(\ref{eq:Neumann-Dirichlet-bc}) follows similarly.
\end{proof-lemma}

\medskip

The stability margin of the platoon can be measured by the real part of $s_1^+$,  the \emph{least stable eigenvalue}.

\medskip

\begin{corollary}\label{cor:s1plus-symmetric}
Consider the eigenvalue problem for the linear
PDE~\eqref{eq:PDE-nominal} with boundary
conditions~(\ref{eq:Dirichlet-bc})
and~(\ref{eq:Neumann-Dirichlet-bc}), corresponding to scenarios~I
and~II respectively. The least stable eigenvalue, denoted by $s_1^+$, satisfies
\begin{align}
s_1^+ & = -\frac{\pi^2k_0}{b_0}\frac{1}{N^2} + O(\frac{1}{N^4}) \quad
      (\text{Dirichlet-Dirichlet})\label{eq:s1plus-DD-nom} \\
s_1^+ & = -\frac{\pi^2
      k_0}{4b_0}\frac{1}{N^2} + O(\frac{1}{N^4}) \quad
      (\text{Neumann-Dirichlet}) \label{eq:s1plus-ND-nom}
\end{align}
as $N\rightarrow\infty$.  \frqed
\end{corollary}
The result shows that the least stable eigenvalue of the closed loop platoon decays as $\frac{1}{N^2}$ with symmetric bidirectional control.

\medskip

\begin{figure}
\psfrag{pde-DD-nom}[l][l][0.9]{PDE~\eqref{eq:PDE-nominal}, D-D}
\psfrag{platoon-DD-nom}[l][l][0.9]{platoon, L-F}
\psfrag{pde-ND-nom}[l][l][0.9]{PDE~\eqref{eq:PDE-nominal}, N-D}
 \psfrag{platoon-ND-nom}[l][l][0.9]{platoon, L}
\psfrag{prediction1}[l][l][0.9]{eq.~\eqref{eq:s1plus-DD-nom}}
\psfrag{prediction2}[l][l][0.9]{eq.~\eqref{eq:s1plus-ND-nom}}
\psfrag{eigenvalues}[l][l][0.9]{$ -Re(s_1^+)$}
\psfrag{N}[l][l][0.75]{$N$}
\begin{center}
\includegraphics[scale=0.45]{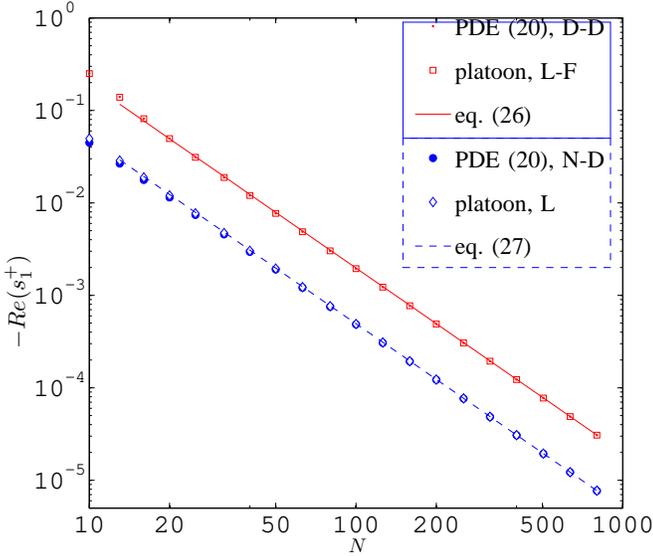}
\end{center}
\caption{\label{fig:eigvaltrendwithn-nominal} Comparison of the least
stable eigenvalue of the closed loop platoon dynamics and that predicted by Corollary~\ref{cor:s1plus-symmetric} with symmetric bidirectional control. There are three plots each for scenarios I and II (corresponding legends are boxed together), and those three should be compared with one another. In the plot legends, ``D-D'' stands for ``Dirichlet-Dirichlet'', ``N-D'' for ``Neumann-Dirichlet'', ``L-F'' for fictitious leader-follower, and ``L'' for fictitious leader. The plot for ``PDE~\eqref{eq:PDE-nominal}, D-D'' should be compared with ``platoon, L-F'' since they both correspond to scenario~I. Similarly, ``PDE~\eqref{eq:PDE-nominal}, N-D'' and  ``platoon, L'' correspond to scenario~II. Note that the predictions~\eqref{eq:s1plus-DD-nom} and~\eqref{eq:s1plus-ND-nom} are valid for $1<<l_c$ (defined in~\eqref{eq:l-critical}), which in this case means for $N>>12$.}
\end{figure}

We now present numerical computations that corroborates this
PDE-based analysis.  Figure~\ref{fig:eigvaltrendwithn-nominal} plots
as a function of $N$ the least stable eigenvalue of the PDE and of the
state-space model of the platoon, as well as the prediction from the
asymptotic formula.  The eigenvalues for the discrete platoon are
obtained by numerically evaluating the eigenvalues of the matrices
$A_{L-F}$ and $A_L$ (see~\eqref{eq:closedloop-platoon-fb}
and~\eqref{eq:closedloop-platoon-f}) with constant control gains
$k^{(f)}_i = k^{(b)}_i = k_0 = 1$ and $b_i = b_0 = 0.5$ for
$i=1,\dots,N$.  The comparison shows that the PDE analysis accurately
predicts the eigenvalue of the state-space model of the platoon
dynamics.

\medskip

%
\begin{figure}
\begin{center}
\psfrag{s+}{$s_l^+$}
\psfrag{s-}{$s_l^-$}
\subfigure[Eigenvalues move toward zero with increasing $N$.]{\includegraphics[scale=0.4, clip = true, trim=0in 0in 0in 0in]{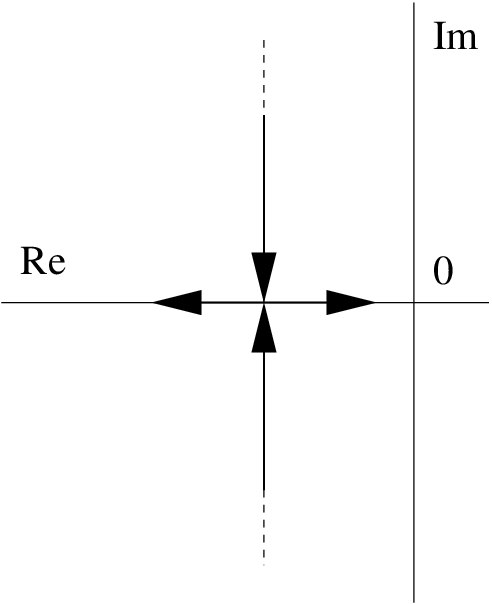}}
\hspace{3 ex}
\subfigure[Mistuning ``exchanges'' stability between $s_l^+$ and $s_l^{-}$.]{\includegraphics[scale=0.4, clip = true, trim=0in 0in 0in 0in]{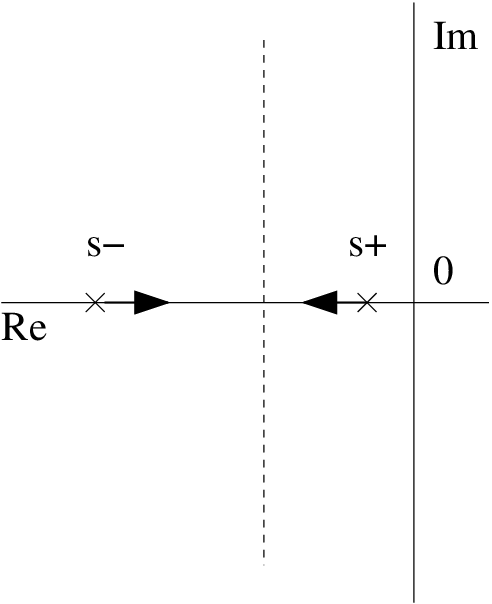}}
\caption{\label{fig:lossofstability-expln}A schematic explaining the loss of stability as $N$ increases and how mistuning ameliorates this loss.}
\end{center}
\end{figure}

Figure~\ref{fig:lossofstability-expln}(a) graphically illustrates the
destabilization by depicting the movement of eigenvalues $s_1^{\pm}$
as $N$ increases.  For sufficiently small values of $N$, the
discriminant $D(1,N)$ is negative and the eigenvalue $s_1^{\pm}$ are
complex. The real part of the eigenvalue depends only on the value
of $b_0$. At a critical value of $N=N_c \eqdef
\frac{\pi\sqrt{2k_0}}{b_0}$, the discriminant becomes zero,
$s_1^+=s_1^-$ and the eigenvalues collide on the real axis. For
values of $N>N_c$ and in particular as $N\rightarrow \infty$, the
eigenvalue $s_1^+$ asymptotes to $0$ while staying real, and $s_1^-$
asymptotes to $-b$.  Their cumulative damping, as reflected in the
sum $s_l^{+} + s_l^{-} = - b_0$, is conserved.  In other words,
$s_1^+$ is destabilized at the expense of $s_1^-$.

\begin{remark} The preceding analysis shows that the loss of
stability experienced with a symmetric bidirectional architecture is
controller independent. The least stable eigenvalue approaches $0$
as $O(1/N^2)$ irrespective of the values of the gains $k_0$ and
$b_0$, as long as they are fixed constants independent of $N$.
Corollary~\ref{cor:s1plus-symmetric} also implies that for the least
stable eigenvalue to be uniformly bounded away from $0$, one has to
increase the control gain $k_0$ as $N^2$. In~\cite{MJ_BB_TAC:05},
the same conclusion was reached for the least stable eigenvalue with
LQR control of a platoon on a circle.  LQR control typically leads
to a centralized architecture, whereas symmetric bidirectional
control is decentralized. It is interesting to note that the least
stable eigenvalue behaves similarly in these distinct architectures.
\frqed
\end{remark}

\comment{
\begin{remark}The result of Lemma~\ref{lem:s1plus-symmetric} also
indicates that the platoon with both a fictitious leader and a
follower (scenario~I) has a higher stability margin than platoon with
only a fictitious leader (scenario~II).  We note, however, that in
scenario~I, the absolute position of the fictitious follow vehicle
must be provided to the last vehicle in the platoon. Therefore
scenario~I requires one extra piece of global information as compared
to the scenario~II.  The result provides a numerical measure of the
benefit of this extra information -- a factor of $4$ improvement in
the closed-loop damping. \frqed
\end{remark}}

\section{Reducing loss of stability by mistuning}\label{sec:mistuning}
In this section, we examine the problem of designing the control
gain functions $k^{(f)}(x), k^{(b)}(x)$ so as to ameliorate the loss
of stability margin with increasing $N$ that was seen in the
previous sections when $k^{(f)}(x) = k^{(b)} \equiv k_0$.
Specifically, we consider the eigenvalue problem for the
PDE~\eqref{eq:PDE-v-1} where the control gains are changed slightly
(mistuned) from their values in the symmetric bidirectional case in
order to minimize the least-stable eigenvalue $s_1^+$.  With
symmetric bidirectional control, one obtains an $O(\frac{1}{N^2})$
estimate for the least stable eigenvalue because the coefficient of
$\frac{\partial^2}{\partial x^2}$ term in PDE~\eqref{eq:PDE-v-1} is
$O(\frac{1}{N^2})$ and the coefficient of $\frac{\partial}{\partial
x}$ term is $0$.  Any asymmetry between the forward and the backward
gains will lead to non-zero $k^{(-)}(x)$ and a presence of
$O(\frac{1}{N})$ term as coefficient of $\frac{\partial}{\partial
x}$.  By a judicious choice of asymmetry, there is thus a potential
to improve the stability margin from $O(\frac{1}{N^2})$ to
$O(\frac{1}{N})$.


We begin by considering the forward and backward position feedback
gain profiles:
\begin{align*}
k^{(f)}(x) & = k_{0} + \epsilon k^{(f,purt)}(x), \\
k^{(b)}(x) & = k_{0} + \epsilon k^{(b,purt)}(x),
\end{align*}
where $\epsilon>0$ is a small parameter signifying the amount of
mistuning and $k^{(f,purt)}(x)$, $k^{(b,purt)}(x)$ are functions
defined over the interval $[0, 2\pi]$ that capture
\emph{perturbation} from the nominal value $k_0$. Define
\begin{align*}
k_s(x) & := k^{(f,purt)}(x)  + k^{(b,purt)}(x), \\
k_m(x) & := k^{(f,purt)}(x) - k^{(b,purt)}(x), 
\end{align*}
so that from~\eqref{eq:kfrontback-x-def},
\begin{align*}
k^{(+)}(x) & = 2k_0 +  \epsilon k_s(x), & & k^{(-)}(x) & = \epsilon k_m(x).
\end{align*}
The mistuned version of the PDE~\eqref{eq:PDE-v-1} is then given by
\begin{equation}\label{eq:linPDE-mistuned}
\frac{\partial^2 v}{\partial t^2}+b_0 \frac{\partial v}{\partial t} =  a_0^2\frac{\partial^2  v}{\partial
x^2}+\epsilon \left[\frac{k_m}{\rho_0} \frac{\partial v}{\partial x}+ \frac{k_s}{2\rho_0^2}\frac{\partial^2  v }{\partial x^2}\right]
\end{equation}
We study the problem of improving the stability margin by judicious
choice of $k_m(x)$ and $k_s(x)$.  The results of our investigation,
carried out in the following sections, provide a systematic
framework for designing control gains in the platoon by introducing
small changes to the symmetric design. \comment{We call this
procedure ``mistuning''.}

\subsection{Mistuning-based design for scenario I}\label{sec:mistuning-perturbation-I}
The control objective is to design mistuning profiles $k_m(x)$ and
$k_s(x)$ to {\em minimize} the least stable eigenvalue $s_1^+$.  To
achieve this, we first obtain an explicit asymptotic formula for the
eigenvalues when a small amount of asymmetry is introduced in the
control gains (i.e., when $\epsilon$ is small).
For scenario~I, the result is presented in the following theorem.
The proof appears in Appendix~\ref{sec:app-proof}.

\begin{theorem}\label{thm:epsilon-mistuning}
Consider the eigenvalue problem for the mistuned
PDE~\eqref{eq:linPDE-mistuned} with Dirichlet boundary
condition~(\ref{eq:Dirichlet-bc}) corresponding to scenario~I.
The $l^\text{th}$ eigenvalue pair is given by the asymptotic formula
\begin{eqnarray*}
s_l^+(\epsilon) & = & \epsilon \frac{l}{2 b_0
N}\int_0^{2\pi}k_m(x)\sin(lx) dx + O(\epsilon^2)+O(\frac{1}{N^2}),\\
s_l^-(\epsilon) & = & - b_0 - \epsilon \frac{l}{2 b_0
N}\int_0^{2\pi}k_m(x)\sin(lx) dx + O(\epsilon^2)+O(\frac{1}{N^2}),
\end{eqnarray*}
that is valid for each $l$ in the limit as $\epsilon\rightarrow 0$
and $N\rightarrow\infty$. \frqed
\end{theorem}






It is apparent from the Theorem above that to minimize the least
stable eigenvalue $s_1^+$, one needs to choose only $k_m$ carefully;
$k_s$ has only $O(\frac{1}{N^2})$ effect. Therefore we choose
$k_s(x) \equiv 0$, or, equivalently, $k^{(f,purt)}(x) = -
k^{(b,purt)}(x)$, which leads to $k_m(x) = 2k^{(f,purt)}(x)$. The
most beneficial control gains are now can be readily obtained from
Theorem~\ref{thm:epsilon-mistuning}, which is summarized in the next
corollary.

\begin{corollary}[Mistuning profile for Scenario~I]\label{cor:mistuning-perturbation}
Consider the problem of minimizing the least-stable eigenvalue of
the PDE~\eqref{eq:linPDE-mistuned} with Dirichlet boundary
condition~(\ref{eq:Dirichlet-bc}) by choosing $k^{(f,purt)}(x)\in
L^\infty([0,2\pi])$ with norm-constraint
$\|k^{(f,purt)}(x)\|_{L^\infty} = \max_{x \in [0  2\pi]}
k^{(f,purt)} (x) = 1$ and $k^{(b,purt)}(x) = - k^{(f,purt)}(x)$.  In
the limit as $\epsilon\rightarrow 0$, the optimal mistuning profile
is given by $k^{(f,purt)}(x) = 2(H(x- \pi)-\frac{1}{2}) $, where
$H(x)$ is the Heaviside function: $H(x) = 1$ for $x\geq 0$ and $H(x)
= 0$ for $x < 0$. With this profile, the least stable eigenvalue is
given by the asymptotic formula
\begin{align*}
s_1^+(\epsilon) = - \frac{4 \epsilon}{b_0N}
\end{align*}
in the limit as $\epsilon\rightarrow 0$ and $N\rightarrow\infty$.
\frqed
\end{corollary}

\medskip

The result shows that even with an \emph{arbitrarily small amount}
of mistuning $\epsilon$, one can improve the closed-loop platoon
damping by a large amount, especially for large values of $N$. The
least-stable eigenvalue $s_1^+$ asymptotes to $0$ as
$O(\frac{1}{N})$ in the mistuned case as opposed to
$O(\frac{1}{N^2})$ in the symmetric case.

\medskip

Figure~\ref{fig:mistuned-platoon-gains}(a) shows the gains for the
individual vehicles (that are obtained from sampling the functions
$k^{(f)}(x)$ and $k^{(b)}(x)$), suggested by
Corollary~\ref{cor:mistuning-perturbation} for a $20$ vehicle
platoon,  with $k_0=1$ and $\epsilon = 0.1$:
\begin{align*}
k^{(f)}_i & = 1 + 0.2(H(\pi-i\delta) - 0.5),\text{ and } \\
k^{(b)}_i & = 1 - 0.2(H(\pi-i\delta) - 0.5),
\end{align*}
where $\delta$ is the desired inter-vehicular spacing in the scaled $y$ coordinates, and is defined in~\eqref{eq:delta-defn}. A confirmation of the
predictions of Corollary~\ref{cor:mistuning-perturbation} is
presented in Figure~\ref{fig:eigvaltrendwithN-mistuning}.
Numerically obtained mistuned and nominal eigenvalues for both the
PDE and the platoon state-space model are shown in the figure, with
mistuned gains chosen as shown in
Figure~\ref{fig:mistuned-platoon-gains}(a). The figure shows that
\begin{enumerate}
\item the platoon eigenvalues match the PDE eigenvalues
accurately over a range of $N$, and
\item the mistuned eigenvalues show large improvement over the nominal
case even though the controller gains differ from their nominal values
only by $\pm 10\%$. The improvement is particularly noticeable for
large values of $N$, while being  significant even for small values of $N$.
\end{enumerate}
For comparison, the figure also depicts the asymptotic eigenvalue
formula given in Corollary~\ref{cor:mistuning-perturbation}.
%

\begin{figure}
\psfrag{front}{$k_i^{(f)}$}
\psfrag{back}{$k_i^{(b)}$}
\psfrag{vehicle index}{vehicle index $i$}
\psfrag{gain}{gain}
\begin{center}
\end{center}
\caption{Mistuned front and back gains $k_i^{(f)}$ and $k_i^{(b)}$ of the vehicles in a platoon with $k_0 = 1$ and $\epsilon = 0.1$. Figure (a) shows the gains chosen according to  Corollary~\ref{cor:mistuning-perturbation} to be optimal for  scenario $I$ for small $\epsilon$:  $k^{(f)}_i = k_0\left(1 +0.1 ( 2H(\pi - i \delta)-1)\right), k^{(b)}_i = k_0\left(1 - 0.1  ( 2H(\pi - i\delta) - 1) \right)$, where $H(\cdot)$ is the Heaviside function and $\delta$ is defined in~\eqref{eq:delta-defn}. Figure (b) shows the optimal mistuned gains for scenario~II with the same parameters, which turns out to be  (see Corollary~\ref{cor:mistuning-perturbation-ND})  $k^{(f)}_i = 1.1k_0$ and $k^{(b)}_i = 0.9 k_0$ for $i=1,\dots,N$. }\label{fig:mistuned-platoon-gains}
\end{figure}

Figure~\ref{fig:lossofstability-expln}(b) graphically illustrates the
mechanism by which mistuning affects the movement of eigenvalues
$s_1^{\pm}$ as $N$ increases. By properly choosing the mistuning
patterns $k_m(x)$ and $k_s(x)$, damping can be ``exchanged'' between
the eigenvalues $s_1^+$ and $s_1^-$ so that the less stable eigenvalue
$s_1^+$ ``gains'' stability at the expense of the more stable
eigenvalue $s_1^-$.  The net amount of damping is preserved, since
$s_1^+ + s_1^- = -b_0$ (as seen from Theorem~\ref{thm:epsilon-mistuning}).

\medskip

\begin{figure}
\psfrag{pde-nom}{symmetric bidi. (PDE)}
 \psfrag{pde-mis}{mistuned (PDE)}
\psfrag{platoon-nom}{symmetric bidi. (platoon)}
\psfrag{platoon-mis}{mistuned (platoon)}
\psfrag{equationxx}{Corollary~\eqref{cor:mistuning-perturbation}}
\psfrag{eigenvalues}{$ -Re(s^+_1)$} \psfrag{N}[][][0.75]{$ N $}
\begin{center}
\includegraphics[scale=0.45]{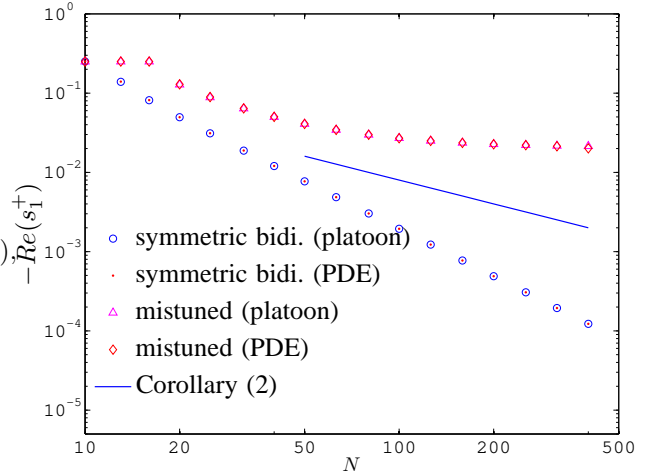}
\end{center}
\caption{\label{fig:eigvaltrendwithN-mistuning} Stability margin improvement by mistuning in Scenario~I. The figure shows the least stable eigenvalue of the closed loop platoon (i.e., of $A_{L-F}$ in
~\eqref{eq:closedloop-platoon-fb}) and of the PDE~\eqref{eq:linPDE-mistuned} with Dirichlet boundary conditions,
with and without mistuning, for a range of values of $N$.  Parameters
for the nominal case are $k_0 = 1$ and $b_0 = 0.5$, and the mistuning amplitude is $\epsilon = 0.1$. The mistuned control gains are shown in Figure~\ref{fig:mistuned-platoon-gains}(a). The legend ``Corollary~\ref{cor:mistuning-perturbation}'' refers to the prediction by Corollary~\ref{cor:mistuning-perturbation} for large $N$.}
\end{figure}

\subsection{Mistuning-based design for scenario II}\label{sec:mistuning-perturbation-II}
For scenario~II, asymptotic formula for the eigenvalue (counterpart
of Theorem~\ref{thm:epsilon-mistuning}) is summarized in the
following theorem.  The proof is entirely analogous to the proof of
Theorem~\ref{thm:epsilon-mistuning}, and is therefore omitted.

\begin{theorem}\label{thm:epsilon-mistunin-ND}
Consider the eigenvalue problem for the mistuned
PDE~\eqref{eq:linPDE-mistuned} with Neumann-Dirichlet boundary
condition~(\ref{eq:Neumann-Dirichlet-bc}) corresponding to
scenario~II. The $l^\text{th}$ eigenvalue pair is given by the
asymptotic formula
\begin{eqnarray*}
s_l^+(\epsilon) & = & -\epsilon \frac{l}{4 b_0 N}\int_0^{2\pi}k_m(x)\sin(\frac{lx}{2}) dx + O(\epsilon^2)+ O(\frac{1}{N^2}),\\
s_l^-(\epsilon) & = & -b_0 + \epsilon \frac{l}{4 b_0
N}\int_0^{2\pi}k_m(x)\sin(\frac{lx}{2}) dx + O (\epsilon^2) + O(\frac{1}{N^2}),
\end{eqnarray*}
that is valid for each $l$ in the limit as $\epsilon\rightarrow 0$
and $N\rightarrow\infty$. \frqed
\end{theorem}

\medskip

As with scenario~I, here again we use the above result to determine
the most beneficial profile $k_m(x)$ for small $\epsilon$:

\medskip

\begin{corollary}[Mistuning profile for Scenario~II]\label{cor:mistuning-perturbation-ND}
Consider the problem of minimizing the least-stable eigenvalue of
the PDE~\eqref{eq:linPDE-mistuned} with Neumann-Dirichlet boundary
conditions~(\ref{eq:Neumann-Dirichlet-bc}) by choosing
$k^{(f,purt)}(x)\in L^\infty([0,2\pi])$ with norm-constraint
$\max_{x \in [0 , 2\pi]} k^{(f,purt)}(x) = 1$, and $k^{(b,purt)}(x)
= - k^{(f,purt)}(x)$.  In the limit as $\epsilon\rightarrow 0$, the
optimal $k^{(f,purt)}$ is given by $k^{(f,purt)}(x) = 1$. With this
profile, the least-stable eigenvalue is given by the asymptotic
formula
\begin{align*}
s_1^+(\epsilon) = - \frac{\epsilon}{b_0N}
\end{align*}
in the limit as $\epsilon\rightarrow 0$ and $N\rightarrow\infty$. \frqed
\end{corollary}

\medskip

\begin{figure}
\psfrag{pde-nom}{symmetric bidi.~(PDE)} \psfrag{pde-mis}{mistuned (PDE)}
\psfrag{platoon-nom}{symmetric bidi.~(platoon)}
\psfrag{platoon-mis}{mistuned (platoon)}
\psfrag{eigenvalues}{$ -Re(s^+_1)$}
\psfrag{N}[][][0.75]{$N$}
\psfrag{equationxx}{Corollary~\eqref{cor:mistuning-perturbation-ND}}
\begin{center}
\includegraphics[scale=0.45]{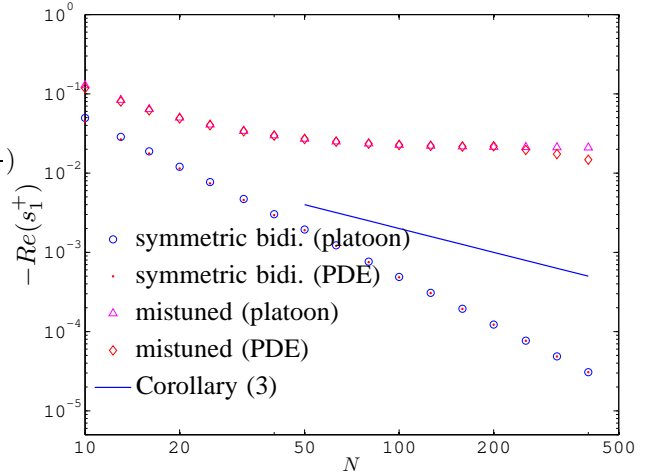}
\end{center}
\caption{\label{fig:eigvaltrendwithN-mistuning-ND}Stability margin improvement by mistuning in scenario~II. The figure shows the least stable eigenvalue of the closed loop platoon (i.e., of $A_L$ in~\eqref{eq:closedloop-platoon-fb}) and of the PDE~\eqref{eq:linPDE-mistuned} with Neumann-Dirichlet b.c., with and
without mistuning, for a range of values of $N$.  The parameters for the nominal case are $k_0 = 1$ and $b_0 = 0.5$, and the mistuning amplitude is $\epsilon = 0.1$. The mistuned control gains that are used are shown in Figure~\ref{fig:mistuned-platoon-gains}(b).
The legend
``Corollary~\ref{cor:mistuning-perturbation-ND}'' refers to the
prediction by Corollary~\ref{cor:mistuning-perturbation-ND} of mistuned
PDE eigenvalues. }
\end{figure}

The result shows that, as in scenario~I, it is possible to improve
the closed-loop stability margin in scenario~II with an arbitrary
small amount of mistuning $\epsilon$ such that the least-stable
eigenvalue $s_1^+$ asymptotes to $0$ as $O(\frac{1}{N})$ in the
mistuned case as opposed to $O(\frac{1}{N^2})$ in the symmetric
case. The gains suggested by Corollary~\ref{cor:mistuning-perturbation-ND}, 
with $k_0=1$ and $\epsilon = 0.1$ are:
\begin{align*}
k^{(f)}_i & = 1.1, & &\text{ and } & k^{(b)}_i & = 0.9,
\end{align*}
which are shown in Figure~\ref{fig:mistuned-platoon-gains}(b).
Numerically obtained least stable eigenvalues for the PDE and
the platoon state-space model for scenario~II are shown in
Fig.~\ref{fig:eigvaltrendwithN-mistuning-ND} for a range of values
of $N$. It is clear from the figure that, as in scenario I, the
mistuned eigenvalues show an order of magnitude improvement over
their values in the symmetric bidirectional case with only $\pm
10\%$ variation.

\medskip

\begin{remark}[Robustness to small changes from the optimal gains]\label{rem:robustness-index}
An advantage of the mistuning design is that  mistuned closed loop eigenvalues are robust to small local discrepancies in the control gains from the optimal ones. This can be seen (for scenario I) from the asymptotic eigenvalue formulas of Theorem~\ref{thm:epsilon-mistuning}, which shows that one would obtain a $O(\frac{1}{N})$ estimate for any choice of $k_m(x)$ such that $\int_0^{2\pi}k_m(x)\sin(x) dx \ne 0$. A similar argument holds for scenario~II.
\end{remark}

\subsection{Simulations}\label{sec:sim}
We now present results of a few simulations that show the
time-domain improvements -- manifested in faster decay of initial
errors -- with the mistuning-based design of control gains.
Simulations were carried out for a platoon of $N=20$ vehicles with
scenario~I, i.e., with fictitious lead and follow vehicles. The
desired gap was $\Delta = 1 $ and desired velocity was $V_d = 5 $.
The initial velocity of every vehicle was chosen as the desired
velocity and the initial position of the $i^\text{th}$ vehicle was
chosen as $Z_i(0) = i\Delta - 0.5$ for $i=\{1,\dots,N\}$. As a
result, the initial relative position error and velocity error of
every vehicle was zero except for the first vehicle, whose relative
position error with respect to the fictitious lead vehicle was
$0.5$.

\medskip

\begin{figure*}
\psfrag{rel-pos-error}[][][0.75]{$Z_{i-1}(t) - Z_{i}(t) - \Delta$}
\psfrag{abs-pos-error}[][][0.75]{$Z_i^d(t) - Z_{i}(t)$}
\psfrag{Time}{$t$}
\begin{center}
\subfigure[Absolute position
errors]{\includegraphics[clip=true,trim=0.0in 0.0in 0.0in 0.0in,scale
= 0.35]{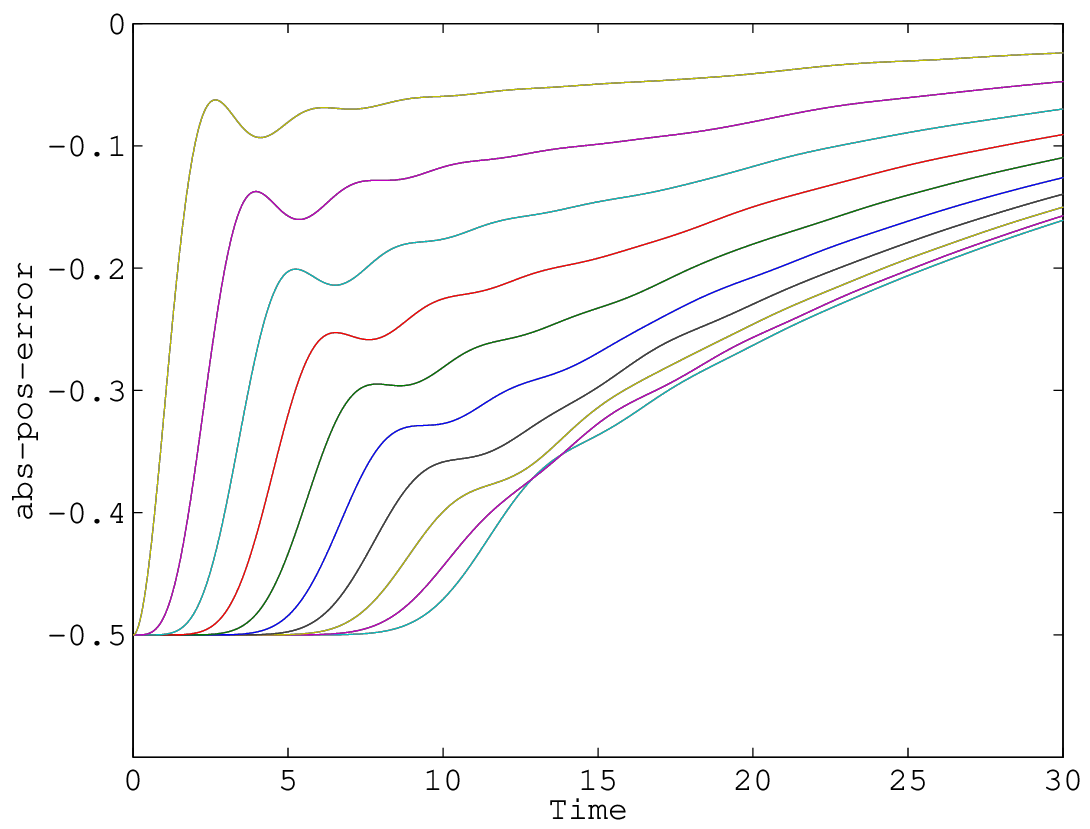}} \hspace{1 ex}
\subfigure[Relative position errors]{\includegraphics[clip=true,trim=0.0in 0.0in 0.0in
0.0in,scale = 0.35]{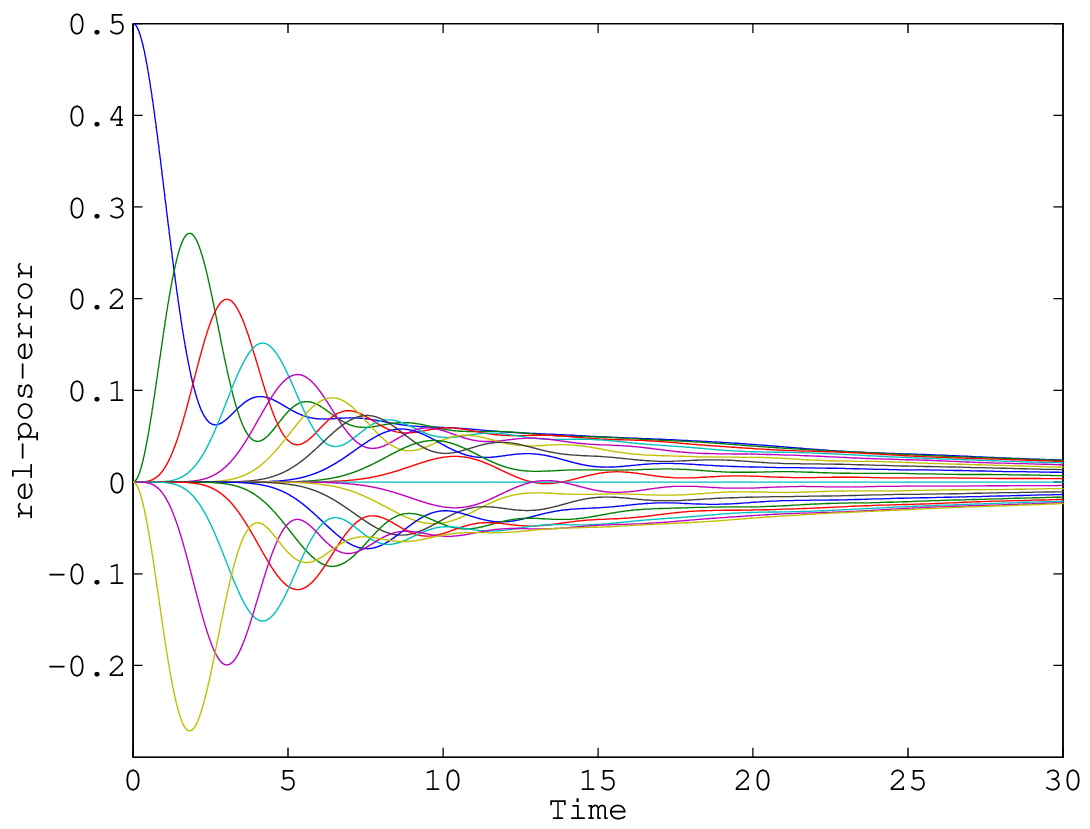}}
\caption{\label{fig:sim-poserror-DD-nom}\emph{Performance of symmetric bidirectional control in time-domain}: time histories of the absolute and relative position errors of the vehicles
in a platoon with symmetric bidirectional control (scenario~I). The control gains are $k^{(f)}_i = k^{(b)}_i = 1$ and $b_i = 0.5$ for every $i=1,\dots,20$. }
\end{center}
\end{figure*}

Figure~\ref{fig:sim-poserror-DD-nom} shows the time-histories of the
absolute and relative position errors of the individual vehicles with
a symmetric bidirectional control, where the control gains were chosen
as $k_i^{(f)} = k_i^{(b)} = 1$ and $b_ i = 0.5$ for
$i=\{1,\dots,20\}$. The absolute position error of the $i^\text{th}$ vehicle is $Z_i - Z_i^d$ and the relative position error is $Z_{i-1} - Z_i -
\Delta$.

\medskip

\begin{figure*}
\psfrag{rel-pos-error}[][][0.75]{$Z_{i-1}(t) - Z_{i}(t) - \Delta$}
\psfrag{abs-pos-error}[][][0.75]{$Z_i^d(t) - Z_{i}(t) $}
\psfrag{Time}{$t$}
\begin{center}
\subfigure[Absolute position errors.]{\includegraphics[clip=true,trim=0.0in 0.0in 0.0in 0.0in,scale = 0.35]{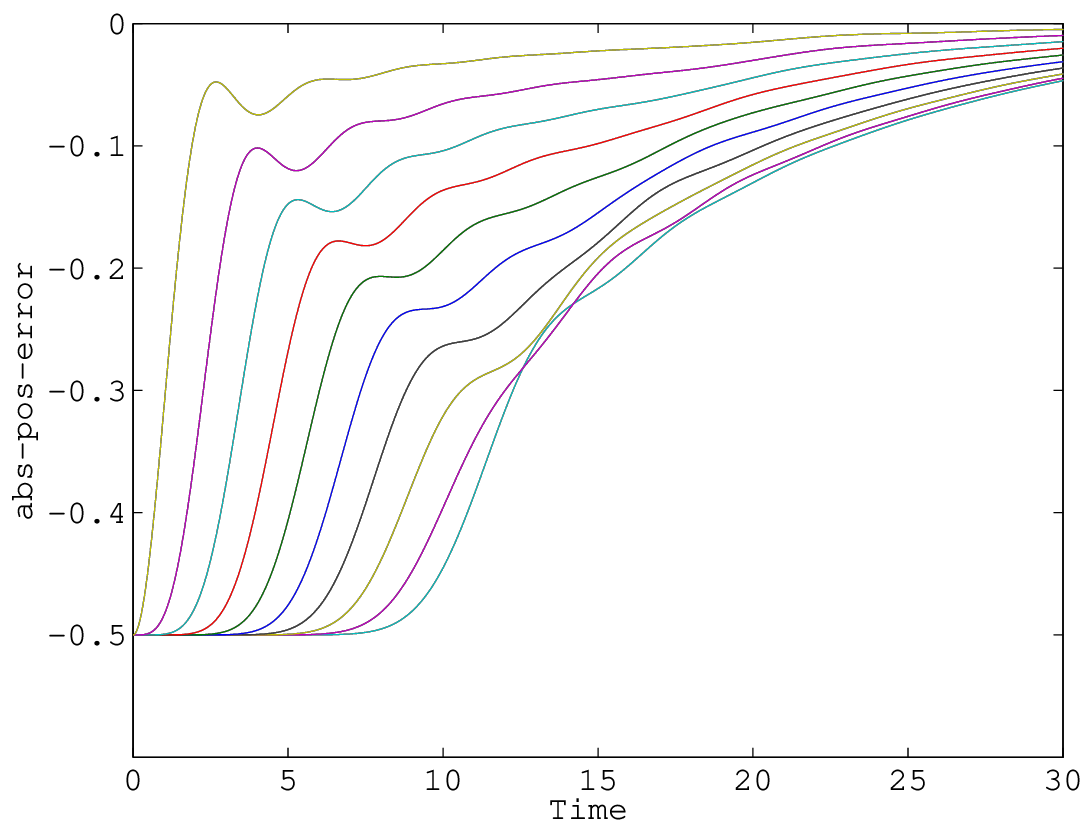}} \hspace{1 ex}
\subfigure[Relative position errors.]{\includegraphics[clip=true,trim=0.0in 0.0in 0.0in 0.0in,scale = 0.35]{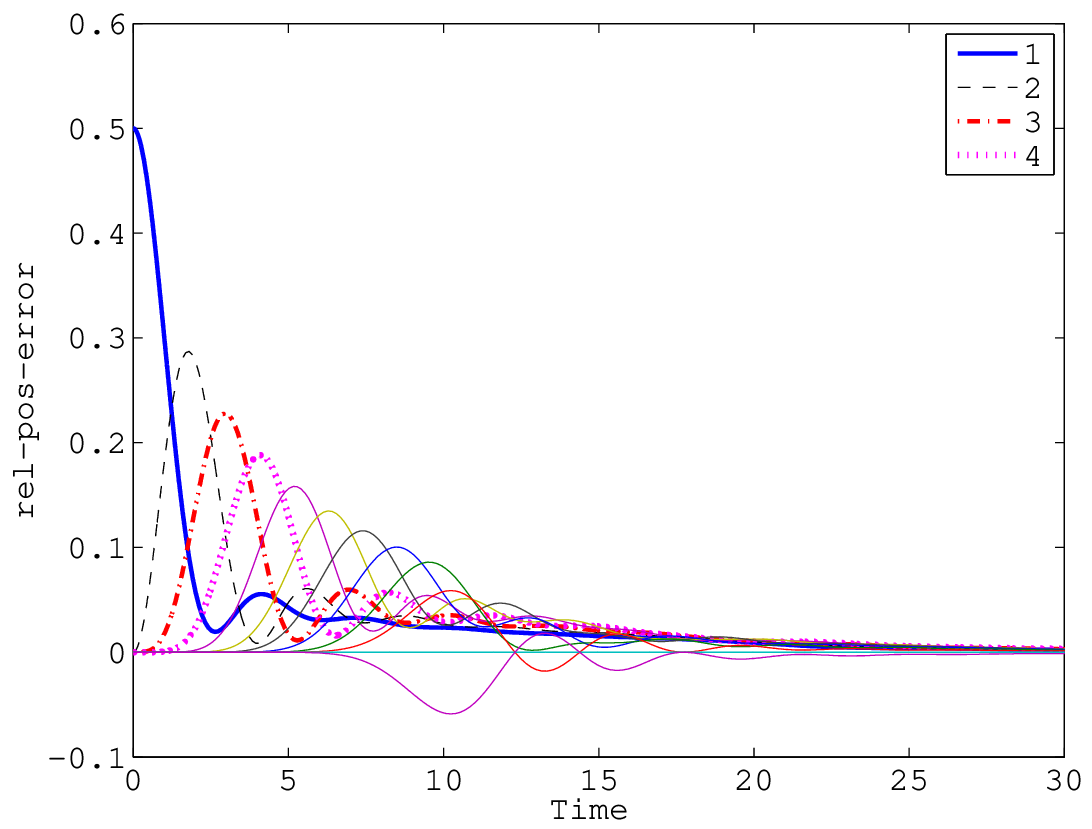}}
\caption{\label{fig:sim-poserror-DD-mis}\emph{Performance of mistuned control 
in time-domain}:
time histories of the absolute and relative position
errors of the vehicles in a platoon (scenario~I) with mistuned
bidirectional control, cf.~Figure~\ref{fig:sim-poserror-DD-nom}. The control gains used are those shown in Figure~\ref{fig:mistuned-platoon-gains}(a). The legends refer to the vehicle indices. }
\end{center}
\end{figure*}

Figure~\ref{fig:sim-poserror-DD-mis} shows the time-histories of
the absolute and relative position errors for the platoon with
mistuned controller gains.  The mistuning gains used for the
simulation are the ones shown in
Figure~\ref{fig:mistuned-platoon-gains}(a) (chosen according to
Corollary~\ref{cor:mistuning-perturbation}) so that maximum and
minimum gains over all vehicles is within $\pm 10\%$ of the nominal
value. On comparing Figures~\ref{fig:sim-poserror-DD-nom}
and~\ref{fig:sim-poserror-DD-mis}, we see that the errors in the
initial conditions are reduced faster in the mistuned case compared
to the nominal case. These observations are consistent with the
improvement in the closed-loop stability margin with the mistuned
design.

\section{Discussion on mistuning design}\label{sec:remarks}
There are several remarks to be made regarding the mistuning based
design. We first comment on the implementation issues, in
particular, on the effect of small platoon size on the proposed
design, and on the information requirements for its implementation.

\medskip

\subsection{Large vs. small $N$}
The PDE model is developed for large $N$.  However, detailed
numerical comparison between the PDE and the discrete state space
model shows that the PDE model provides quantitatively correct
predictions even for small values of $N$ (see
Figures~\ref{fig:eigvaltrendwithn-nominal},~\ref{fig:eigvaltrendwithN-mistuning}
and~\ref{fig:eigvaltrendwithN-mistuning-ND}).  The PDE has an
infinite number of eigenvalues as opposed to a finite number for the
discrete platoon.  So, one can not expect an exact match.  However,
PDE eigenvalues exactly match the least stable and other dominant
eigenvalues of the discrete platoon (see
Figure~\ref{fig:PDEvalidation-spectra1} and Figure~\ref{fig:sixeigenvalues_mistuning_01}). In a similar vein, the benefits
of mistuning are also realized for small values of $N$.
For example, when the number of vehicles is $20$, a mistuning of $\pm 10\%$ results in
an improvement in the stability margin -- as measured by the real part of the least stable eigenvalue -- of $150\%$ (from $-0.0491$ to $-0.1281$ ) in scenario
I and an improvement of $400\%$ (from $-0.012$ to $-0.05$) in
scenario~II over the symmetric case.

\subsection{Information requirements}\label{sec:implementation}
In order to implement the beneficial mistuned controller gains
designed above, every vehicle needs the following information (in
addition to what is needed to use a symmetric bidirectional
control): (1) the mistuning amplitude $\epsilon$, and (2) in
scenario~I, whether it is in the front half of the platoon or not.
This information can be provided to the vehicles in advance. In
scenario~II, only the value of $\epsilon$ is needed.

It is possible that due to vehicles leaving and joining the platoon,
information on whether a  vehicle belongs to the front half of the
platoon may become erroneous with time, especially for the vehicles that are
close to the middle. In scenario~I, such error may lead to a
non-optimal gains used by the vehicles. However, since the
improvement in closed loop stability margin due to mistuning is
robust to small deviations in the gains from the optimal ones (see
Remark~\ref{rem:robustness-index}), errors in determining whether a
vehicle belongs to the front half of the platoon or not will not
greatly affect the improvement in stability margin. Note that in
scenario~II this issue does not even arise.

\subsection{Large asymmetry}\label{sec:large-asymmetry}
\begin{figure}
\psfrag{eigenvalues}[][][0.75]{$ - Re(s_l^+), l = 1,2,\dots,6$}
\psfrag{N}[][][0.75]{$N$}
\psfrag{pde}[][][0.6]{PDE}
\begin{center}
\includegraphics[scale=0.4]{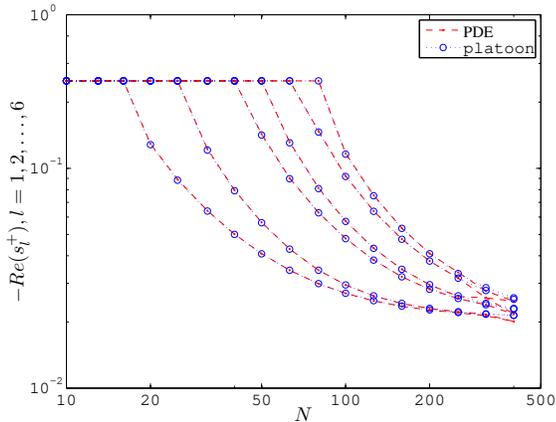}
\end{center}
\caption{\label{fig:sixeigenvalues_mistuning_01}The real parts of six eigenvalues (closest to $0$) of the closed
loop platoon dynamics for Scenario I, and their comparison with the
PDE eigenvalues with Dirichlet-Dirichlet boundary conditions, with
controller gains mistuned  as those shown in Figure~\ref{fig:mistuned-platoon-gains}.  As predicted by the S-L theory, the least stable eigenvalue stays the least stable, although
eigenvalues that are more stable merge with it as $N$ increases.}
\end{figure}

Although the mistuning profiles described in
Corollaries~\ref{cor:mistuning-perturbation}
and~\ref{cor:mistuning-perturbation-ND} are optimal in the limit as
$\epsilon \to 0$, one would like to be able to use them with
somewhat larger values of $\epsilon$ to realize the benefit of
mistuning. To do so, one has to preclude the possibility of
``eigenvalue cross-over'', i.e., of the second ($s_2^+$) or some
other marginally stable eigenvalue from becoming the least stable
eigenvalue in the presence of mistuning.  It turns out that such a
cross-over is ruled out as a consequence of the Strum-Liouville
(S-L) theory for the elliptic boundary value problems.  The standard
argument relies on the positivity of the eigenfunction corresponding
to $s_1^+$; the reader is referred to~\cite{Evans:98} for the
details. Figure~\ref{fig:sixeigenvalues_mistuning_01} verifies this
numerically by depicting the six eigenvalues closest to $0$ (for
both the PDE and the discrete platoon) as a function of $N$ when
mistuning is applied.

\medskip

\subsection{Sensitivity to disturbance}\label{sec:disturbance}
Automated platoons suffer from high sensitivity to external
disturbances; which is referred to as ``string instability'' or
``slinky-type
effects''~\cite{SwaroopHedrick_stringstability_TAC:96,DS_KH_CC_PI_VSD:94,Zhang_IntelligentCruise_TVT:99}. Here we provide numerical evidence that mistuning also helps in reducing the sensitivity to disturbances.

When external disturbances are present, we model the dynamics of vehicle $i$ by $\ddot{Z}_i = U_i +  W_i$, where $W_i$ is the external disturbance acting on the vehicle.  In the $y$ coordinates, the vehicle dynamics become $\ddot{\tilde{y}}_i = u_i + w_i$, where $w_i \eqdef 2\pi W_i/L $. In scenario~I, the state space model of the entire platoon becomes,
\begin{align}\label{eq:ABCD}
\dot{\mbf \psi} & = A_{L-F}\;{\mbf \psi} + \underbrace{\matt{\mbf{0} \\ I}}_{\mathcal{B}}\mbf{w}, &        \mbf{e} & = C \mbf{\psi}
\end{align}
where $\mbf{\psi} = [{\mbf{\tilde y}}^T, {\mbf{\tilde v}}^T]^T$,
$A_{L-F}$, $\mbf{w} =[w_1,w_2,\dots,w_N]^T$, and $\mbf{e} \eqdef
[e_1^{(f)},\dots,e_N^{(f)}]^T$ is a vector of front spacing errors $e_i^{(f)} \eqdef \tilde{y}_{i-1} - \tilde{y}_i$.

\begin{figure}
\psfrag{symm}{symmetric}
\psfrag{mistuned}{mistuned}
\psfrag{H}{$\|G_{we}\|_\infty$}
\psfrag{N}{$N$}
\begin{center}
\includegraphics[scale=0.4]{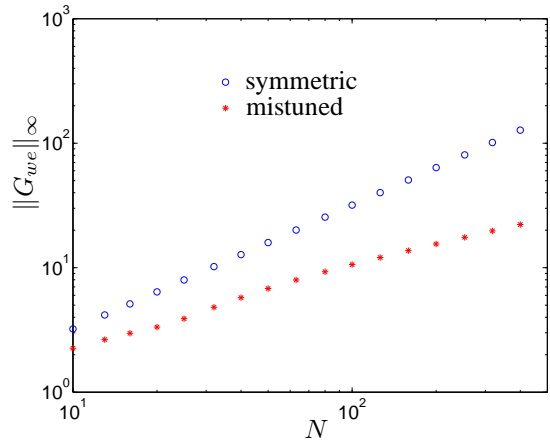}
\caption{\label{fig:H-infinity}$H_\infty$ norm of the transfer
function $G_{we}$ from disturbance $\mbf w$ to spacing error $\mbf
{e}$ in~\eqref{eq:ABCD}, with and without mistuning, for scenario~I.  The mistuned gains used are shown in Figure~\ref{fig:mistuned-platoon-gains}(a). Norms are computed using the Control Systems Toolbox in MATLAB$^\text{\copyright}$.}
\end{center}
\end{figure}

The $H_\infty$ norm of the transfer function $G_{we}$ from the
disturbance $\mbf{w}$ to the inter-vehicle spacing errors $\mbf{e}$
is a measure of the closed loop's sensitivity to external
disturbances~\cite{Seiler_disturb_vehicle_TAC:04,PB_JPH_CDC:05}.
Figure~\ref{fig:H-infinity} shows a plot of the $H_\infty$ norm of
$G_{we}$ as a function of $N$, with and without mistuning. The
mistuning profile used is the same as the one used for the
eigenvalue trends reported in
Figure~\ref{fig:eigvaltrendwithN-mistuning}. It is clear from the
figure that $\pm 10\%$ mistuning results in large reduction of the
$H_\infty$ norm of $G_{we}$. Although this reduction is more
pronounced for large $N$, it is still significant for small $N$. In
particular, for $N=20$, a $10\%$ mistuning yields approximately
$50\%$ reduction in the $H_\infty$ norm (from 6.69 to 3.38).

Apart from the $H_\infty$ norm of $G_{de}$, there are other ways to measure sensitivity to disturbances. In~\cite{RM_JB_TACsubmit:08}, the transfer function from disturbance acting on the lead vehicle to spacing error on the $i^\text{th}$ vehicle is analyzed. Detailed analysis of the effect of mistuning on sensitivity to disturbances will be a subject of future work.

\section{Conclusion}
We developed a PDE model that describes the closed loop dynamics of an $N$-vehicle platoon
 with a decentralized bidirectional control architecture.  Analysis
of the PDE model revealed several important features of the problem.
First, we showed that when every vehicle uses the same controller
with constant gain that is independent of $N$ (the so-called
symmetric bidirectional architecture), the least stable eigenvalue
of the closed loop decays to $0$ as $O(\frac{1}{N^2})$.  Second, and
more significantly, analysis of the PDE suggested a way to
ameliorate the progressive loss of stability with increasing $N$, by
introducing small amounts of ``mistuning'', i.e., by changing the
controller gains from their nominal symmetric values. We proved that
with arbitrary small amounts of mistuning, the decay of the least
stable closed loop eigenvalue can be improved to ${O}(\frac{1}{N})$.
Several comparisons with the numerically computed eigenvalues of
state-space model of the platoon confirm the predictions of the
PDE-based analysis.

\medskip

Although the PDE model is derived under the assumption that the
number of vehicles, $N$, is large, in practice the PDE provides
quantitatively correct predictions for the discrete platoon dynamics
even for relatively small values of $N$.
The amount of information that is needed to implement the mistuned
control gains (over that in the symmetric bidirectional
architecture) is quite small and need to be provided only once.
Furthermore, the stability improvement due to mistuning is robust to
small errors (between the actual gains used and the optimal mistuned
gains) that may occur in practice due to changes in the number of
vehicles in the platoon over time.

\medskip

The advantage of the PDE formulation is reflected in the ease with which
the closed loop eigenvalues are obtained for two different boundary
conditions, with lead and follow vehicles as well as with only a
lead vehicle. Certain important aspects of the problem, such as the beneficial nature of
forward-backward asymmetry in control gains, is revealed by the PDE
while they are difficult to see with the (spatially) discrete,
state-space model.

Numerical calculations show that the mistuning
design also reduces sensitivity to disturbances of the closed-loop platoon.
Analysis of the beneficial effect of mistuning in reducing sensitivity to external disturbances is a subject of future research. In the future, we also plan to examine PDE-based models for modeling and analysis of fleet of vehicles as in $2$ or $3$ spatial dimensions.

\bibliography{vehicular_platoon}
\ifthenelse{\equal{\PaperORReport}{paper}}{
\begin{biography}[{\includegraphics[width=1in,height=1.25in,clip,keepaspectratio
]{./epsfiles/Baroo.eps}}]{Prabir Barooah}
Prabir Barooah is an Assistant Professor of Mechanical and Aerospace Engineering at the University of Florida, where he has been since October 2007. He was born in Jorhat, Assam (India) in 1975. He received his B.Tech and M.S. degrees in Mechanical Engineering from the Indian Institute of Technology, Kanpur in 1996 and University of Delaware, Newark, in 1999, respectively. From 1999 to 2002 he worked at United Technologies Research Center as an associate research engineer. He received the Ph.D. degree in Electrical and Computer Engineering from the University of California, Santa Barbara in 2007.  His research interests include system identification, decentralized and cooperative control of multi-agent systems, and estimation in large-scale sensor networks. Dr. Barooah is the winner of the best paper award at the 2nd Int.~Conf.~on Intelligent Sensing and Information Processing, and a 2003 NASA group achievement award.
\end{biography}
\begin{biography}[{\includegraphics[width=1in,height=1.25in,clip,keepaspectratio
]{./epsfiles/Mehta.eps}}]{Prashant G. Mehta}
Dr. Prashant G. Mehta is an Assistant Professor at the Department of Mechanical \& Industrial Engineering, University of Illinois at Urbana-Champaign.  He received his Ph.D. in Applied Mathematics from Cornell University in 2004.  Prior to joining UIUC, he was a research engineer at the United Technologies Research Center (UTRC).  At UTRC, he was recognized with an outstanding achievement award for his contributions in developing Dynamical Systems methods to obtain practical solutions to problems in aero-engines.  His research interests are in Dynamical Systems and Control including its applications to stochastic dynamics and control of network problems, fundamental limitations in control of non-equilibrium dynamic behavior, and multi-scale and symmetry-aided analysis of interconnected systems.
\end{biography}
\begin{biography}[{\includegraphics[width=1in,height=1.25in,clip,keepaspectratio
]{./epsfiles/Hespan.eps}}]{Jo\~{a}o P. Hespanha}
João P. Hespanha received the Licenciatura in electrical and computer engineering from the Instituto Superior T\'{e}cnico, Lisbon, Portugal in 1991 and the M.S. and Ph.D. degrees in electrical engineering and applied science from Yale University, New Haven, Connecticut in 1994 and 1998, respectively. He currently holds a Professor position with the Department of Electrical and Computer Engineering, the University of California, Santa Barbara. From 1999 to 2001, he was an Assistant Professor at the University of Southern California, Los Angeles. Dr. Hespanha is the associate director for the Center for Control, Dynamical-systems, and Computation (CCDC) and an executive committee member for the Institute for Collaborative Biotechnologies (ICB), an Army sponsored University Affiliated Research Center (UARC). His research interests include hybrid and switched systems; the modeling and control of communication networks; distributed control over communication networks (also known as networked control systems); the use of vision in feedback control; stochastic modeling in biology; the control of haptic devices, and game theory. He is the author of over one hundred technical papers and the PI and co-PI in several federally funded projects. Dr. Hespanha is the recipient of the Yale University’s Henry Prentiss Becton Graduate Prize for exceptional achievement in research in Engineering and Applied Science, a National Science Foundation CAREER Award, the 2005 Automatica Theory/Methodology best paper prize, and the best paper award at the 2nd Int. Conf. on Intelligent Sensing and Information Processing. Since 2003, he has been an Associate Editor of the IEEE Transactions on Automatic Control.
\end{biography}
}
{
}


\appendices



\section{Technical results}\label{sec:appendix}

\subsection{Solution properties of PDE~\eqref{eq:PDE-v-1}.}\label{sec:pde-properties}
In this section, we use the semigroup theory to obtain results on
well-posedness of the PDE~\eqref{eq:PDE-v-1}.  To apply these
methods, we first re-write the PDE as a first order evolution
equation:
\begin{equation}\label{eq:eq_evolu_0}
\begin{array}{ccc}
\frac{\partial \tilde\rho}{\partial t} &=& -\rho_0 \frac{\partial
 v}{\partial x} \\ \frac{\partial  v}{\partial t} &=&
-\left[ \frac{1}{\rho_02}k_1(x) \tilde\rho +
   \frac{1}{2\rho_03}\frac{\partial }{\partial x}(\tilde \rho
   k_0(x)) +  b v \right]
   \end{array} \eqdef A \left[ \begin{array}{c} \tilde{\rho} \\ v
   \end{array} \right],
\end{equation}
where $A$ is a linear operator; $k_0(x)\eqdef k^{+}(x)$ and
$k_1(x)\eqdef k^{-}(x) - \frac{1}{2\rho_0} \frac{d k^{+}}{dx}(x)$.
We will assume these coefficients $k_0(x), k_1(x) \in
L^\infty([0,2\pi])$ and $k_0(x)>0$.  $\tilde\rho$ has the units of
and the physical interpretation of density perturbation.

Using~\eqref{eq:eq_evolu_0}, we denote the initial/boundary value
problem as:
\begin{eqnarray}
\dot{z}(x,t) &=& Az(x,t) \quad \textrm{for} \;\; x\in X, \;\; t>0\nonumber\\
z(x,0) &=& z_0(x),\label{eq:eq_evolu}
\end{eqnarray}
where $z(x,t)\eqdef [\tilde{\rho}(x,t),v(x,t)]$, $z_0(x) =
[\tilde{\rho}_0(x),v_0(x)]$ and $A$ is defined
in~\eqref{eq:eq_evolu_0}; $\tilde{\rho}_0$ and $v_0$ will be
assumed to functions in appropriately defined Banach spaces.  The
main goal of this section will be to show that the solution for the
linear problem~\eqref{eq:eq_evolu_0} can be expressed in terms of a
$C^0$ semigroup provided eigenvalues of the operator $A$ satisfy
appropriate bounds.
We begin with a discussion of the notation.

{\em Preliminaries and Notation.}  We denote $z\eqdef
[\tilde{\rho},v]$, $L^2(X)$ denotes the Hilbert space of
square integrable functions on $X$ ($\|v\|_{L^2}^2 \eqdef
\int v^2 dx$), $H^k$ denotes the Sobolev space of functions
such that derivatives up to $k^{th}$-order exist in a weak sense and
belong to $L^2(X)$ (the Sobolev norm is denoted by
$\|\cdot\|_{H^k}$), and $H_0^1$ denotes the Sobolev space $H^1$ of
functions that satisfy the Dirichlet boundary condition. We denote
$Z\eqdef L^2 \times L^2$, and equip it with a norm $\|\cdot\|$.
Let ${\cal D}(A) \eqdef H^1\times(H^1_0\cap L^2)$ and consider the
right hand side of evolution equation~\eqref{eq:eq_evolu_0} as an
unbounded but closed densely defined linear operator
\begin{equation}
A:{\cal D}(A)\subset Z \rightarrow Z.
\end{equation}
A real number $s$ belongs to ${\rho}(A)$, the resolvent set for $A$,
provided the operator $sI-A:{\cal D}(A) \rightarrow Z$ is 1-1 and
onto.  For $s\in\rho(A)$, the resolvent operator
$R_s\eqdef(sI-A)^{-1}$. Finally, we recall that a  one-parameter
family of linear operators $\{S(t)\}_{t\ge 0}$ is a $C^0$-semigroup
if 1) $S(0)z=z$ for all $z\in Z$, 2) $S(t+s)z=S(t)S(s)z$ for all
$t,s\ge 0$ and $z\in Z$, and 3) the mapping $t\rightarrow S(t)z$ is
continuous from $[0,\infty)$ into $Z$.  A $C^0$ semigroup is a
contraction semigroup if $\|S(t)z\|\le  \|z\|$ for all $t\ge 0$. The
Hille-Yosida theorem states that a closed densely defined linear
operator $A$ is the generator of a
contraction semigroup if and only if
\begin{equation}\label{eq:resolvent_bound}
(0,\infty)\subset \rho(A)\quad\textrm{and}\quad \|R_s z\|\le
\frac{1}{s} \|z\| \quad \forall z\in Z.
\end{equation}

Our strategy will be to apply Hille-Yosida theorem to deduce
solution properties of the evolution equation~\eqref{eq:eq_evolu}.
Following closely the development in~\cite{Evans:98}, there are
three steps to accomplish this: 1) we show that $A$ is a densely
defined closed linear operator on $Z$, 2) characterize the resolvent
set by considering the eigenvalue problem, and 3) show the
bound~\eqref{eq:resolvent_bound} for the resolvent.  Step 2 will
lead to an eigenvalue problem, whose analysis and optimization is
the subject of this paper. We present details for the three steps
next:
\begin{enumerate}
\item The domain of $A$, ${\cal D}(A)$, is dense in $Z$ because
$H^1$ is dense in $L^2$.  To show $A$ is closed, consider a sequence
$\{\tilde{\rho}_m,v_m\}\subset {\cal D}(A)$ such that
\begin{eqnarray}
(\tilde{\rho}_m,v_m) &\stackrel{Z}{\rightarrow}& (\tilde{\rho},v)\label{eq:closed_cond_1}\\
A(\tilde{\rho}_m,v_m) &\stackrel{Z}{\rightarrow}&
(f,g),\label{eq:closed_cond_2}
\end{eqnarray}
where the arrow notation denotes the fact that the convergence is in
$Z=L^2\times L^2$.  Since $v_m \stackrel{L^2}{\rightarrow}
v$ so $-{\rho}_0\frac{\partial v}{\partial x}=f\in
L^2$, i.e., $v\in H^1$. Now, $\{v_m\}$ is Cauchy in
$L^2$ by~\eqref{eq:closed_cond_1} and $\{\frac{\partial
 v_m}{\partial x}\}$ is Cauchy in $L^2$
by~\eqref{eq:closed_cond_2} and
\begin{equation}
\|v_m- v_l\|_{H^1}\le C\left(\|\frac{\partial
v_m}{\partial x}-\frac{\partial v_l}{\partial
x}\|_{L^2} + \|v_m-v_l\|_{L^2} \right),
\end{equation}
so $\{v_m\}$ is Cauchy in $H^1$ and $v_m
\stackrel{H^1}{\rightarrow} v$.  By repeating essentially
the same argument, one also finds that $\tilde{\rho}\in H^1$ and
$\tilde{\rho}_m \stackrel{H^1}{\rightarrow} \tilde{\rho}$.
Consequently, $A(\tilde{\rho}_m,v_m)
\stackrel{Z}{\rightarrow} A(\tilde{\rho},v)$ and
$A(\tilde{\rho},v)=(f,g)$.
\item Let $s>0$, $(f,g)\in Z= L^2 \times L^2$, and consider the
operator equation
\begin{equation}\label{eq:resolvent_eq1}
(sI - A)\left[ \begin{array}{c} \tilde{\rho} \\ v
\end{array}\right]=\left[
\begin{array}{c} f\\g \end{array}\right].
\end{equation}
This is equivalent to two scalar equations
\begin{align}
s \tilde{\rho} + \rho_0 \frac{\partial v}{\partial x} = f\quad (\tilde{\rho}\in
L^2\cap H^1),\label{eq:eq_rs_0}\\
s v + \left[ \frac{1}{\rho_0^2}k_1(x) \tilde{\rho} +
   \frac{1}{2\rho_0^3}\frac{\partial }{\partial x}(
   k_0(x)\tilde{\rho}) +  b v \right]= g \quad (v \in
L^2\cap H^1_0).\label{eq:eq_rs_1}
\end{align}
Using the first equation to write $s\tilde{\rho} = -\rho_0
\frac{\partial v}{\partial x} + f$, this implies
\begin{equation}\label{eq:resolvent_eq1_1}
s^2 v+ b s v + L v = h,
\end{equation}
where
\begin{equation}\label{eq:resolvent_eq2}
L v \eqdef \frac{1}{2\rho_0^2}\frac{\partial}{\partial
x}(-k_0(x) \frac{\partial  v}{\partial x}) -
\frac{1}{\rho_0}(k_1(x) \frac{\partial  v }{\partial x})
\end{equation}
is an elliptic operator (because $k_0 (x)>0$ for all $x\in X$) and
$h=sg - \frac{1}{2\rho_0^3}\frac{\partial }{\partial x}(
   k_0(x)f)-\frac{1}{\rho_0^2}k_1(x)f$ (note that $h \in
   H^{-1}(X)$).
Consequently, solutions of~(\eqref{eq:resolvent_eq1}) can be studied
in terms of solutions of~(\eqref{eq:resolvent_eq2}). The spectrum of
$A$ is completely characterized by the spectrum of $L$. We will
obtain spectral bounds, dependent upon $k_0(x)$ and $k_1(x)$, in the
following sections.  In particular, we will establish that
$\textrm{Real}[s]<\alpha$ for some $\alpha<0$ and thus
$\rho(A)\supset (\alpha,\infty)$. For $k_1(x)=0$, its turns out that
$[0,\infty)\subset\rho(A)$ for any choice of positive $k_0(x)$ (this
is also clear from the symmetric eigenvalue
problem~\eqref{eq:resolvent_eq2}).
\item
If a positive $s\in\rho(A)$, there exists a unique solution
$(\tilde{\rho},v)\in Z$
for~\eqref{eq:eq_rs_0}-\eqref{eq:eq_rs_1} via the theory of elliptic
operators: solve~\eqref{eq:resolvent_eq1_1} to obtain $v\in
H_0^1$ and $s\tilde{\rho} = -\rho_0 \frac{\partial
 v}{\partial x} + f$. We write the solution as
$(\tilde{\rho},v)=R_s(f,g)$, define a bilinear form
\begin{equation}
B[\tilde{\rho},s]\eqdef \frac{1}{2\rho_0^4}\int_X
k_0(x)\tilde{\rho}(x)s(x)dx,
\end{equation}
for $\tilde{\rho},s\in L^2$ and consider an equivalent norm (on $Z$)
for solutions $(\tilde{\rho}, v)$ as:
\begin{equation}
\|(\tilde{\rho}, v )\|\eqdef B[\tilde{\rho},\tilde{\rho}] +
\| v \|_{L^2}
\end{equation}
To obtain the resolvent bound, we multiply~\eqref{eq:eq_rs_1} by
$ v $ and use integration by parts:
\begin{multline*}
s(\| v \|_{L^2}+B[\tilde{\rho},\tilde{\rho}]) +
b\|v \|_{L^2}+\frac{1}{\rho^2}\int k_1(x)\tilde{\rho}
v dx \\
= \int  g v dx + B[\tilde{\rho},f].
\end{multline*}
In general, the bound depends upon $k_1(x)$.  For $k_1(x)=0$, we
have
\begin{multline*}
s\|(\tilde{\rho},v)\|^2\le
(s+b)\| v\|_{L^2}+B[\tilde{\rho},\tilde{\rho}]) = \\
\int g v dx + B[\tilde{\rho},f] \le
\|(f,g)\|\cdot\|(\tilde{\rho},v)\|,
\end{multline*}
where the first inequality holds because $s>0$ and $b>0$ and the
last inequality follows from the generalized Cauchy-Schwarz
inequality.  As a result, $\|R_s(f,g)\|\le\frac{1}{s}\|(f,g)\|$ and
$\|R_s\|\le\frac{1}{s}$.

For the general case where $k_1(x)$ is not identically zero, one
expresses the operator
\begin{equation}
A=A_0 + \tilde{A},
\end{equation}
where
\begin{align*}
A_0 \left[ \begin{array}{c} {\tilde{\rho}} \\ {\tilde{v}}
\end{array} \right]  & =
\left[ \begin{array}{cc} 0 & -\rho_0 \frac{\partial}{\partial x} \\
-\frac{1}{2\rho_0^3}\frac{\partial }{\partial x}(
   k_0(x)\cdot) & -b
\end{array} \right] \left[ \begin{array}{c} {\tilde{\rho}} \\ {v} \end{array} \right], \\
\tilde{A} \left[ \begin{array}{c} {\tilde{\rho}} \\
{v} \end{array} \right] & =
\left[ \begin{array}{cc} 0 &  0 \\
-\frac{1}{\rho_0^2}k_1(x) & 0
\end{array} \right] \left[ \begin{array}{c} {\tilde{\rho}} \\ {v} \end{array}
\right].
\end{align*}
In words, $A_0$ is the operator with $k_1(x)=0$ and $\tilde A$ is
the operator due to $k_1(x)$.  We note that $\tilde A$ is a bounded
perturbation of $A_0$ (on $Z$).  We have already showed the
existence of a $C^0$-semigroup for $A_0$.  For the general operator
$A$, the existence follows from using a perturbation theorem
(see Theorem~1.1 in Ch.~3 of~\cite{MR710486}).
\end{enumerate}

\subsection{Proof of Theorem~\ref{thm:epsilon-mistuning}}\label{sec:app-proof}
\begin{proof-theorem}{\ref{thm:epsilon-mistuning}}
The spatial inhomogeneity introduced by the $x$-dependent coefficients
$k_m(x)$ and $k_s(x)$ destroy the spatial invariance of the nominal
PDE~\eqref{eq:PDE-nominal}.  Hence, the Fourier basis --
eigenfunctions of the Laplacian -- no longer lead to a diagonalization
of the mistuned PDE.  The methods of
section~\ref{sec:instability-analysis} thus need to be suitably
modified.  In order to compute the eigenvalues for the mistuned
PDE~\eqref{eq:linPDE-mistuned}, we take a Laplace transform
of~\eqref{eq:linPDE-mistuned} and get
\begin{equation}\label{eq:PDE_mis1}
-a_0^2\frac{\partial^2 \eta}{\partial x^2} + s^2\eta + b_0 s\eta =
\epsilon\left[\frac{k_m}{\rho_0} \frac{\partial \eta}{\partial x}+
\frac{k_s}{2\rho_0^2}\frac{\partial^2\eta}{\partial x^2} \right],
\end{equation}
where $\eta(x)$ is the Laplace transform (with respect to $t$) of
$\tilde{v}(x,t)$. We are interested in eigenvalues
of~\eqref{eq:PDE_mis1} with Dirichlet boundary conditions, i.e., the
values of $s$ for which a solution to the homogeneous
PDE~\eqref{eq:PDE_mis1} exists with boundary conditions
$\eta(0)=\eta(2\pi)=0$.  To obtain these eigenvalues, we use a
perturbation method expressing the eigenfunction and eigenvalue in a
series form:
\begin{align}
\eta(x) & =  \eta_0(x)+\epsilon \eta_1(x)+ O(\epsilon^2), &  s  & = r_0 + \epsilon r_1+ O(\epsilon^2). \label{eq:perturb}
\end{align}
We note that $\epsilon \, r_1$ denotes the perturbation to the nominal eigenvalue
$r_0$ as a result of the mistuning.  Substituting~\eqref{eq:perturb}
in~\eqref{eq:PDE_mis1} and doing an $O(1)$ balance, we get
\begin{equation}
O(1):\;\;\;\; -a_0^2 (\eta_0)_{xx} + r_0^2 \eta_0 + b r_0 \eta_0 =
0,
\end{equation}
whose eigen-solution is given by
\begin{align*}
\eta_0 & =  d_l \sin(\frac{l x}{2}), & r_0 &= s_l^{\pm}(0),
\end{align*}
where $l=1,2,\hdots$,  $d_l$ is an arbitrary real constant, and
$s_l^{\pm}(0)$ is given by~\eqref{eq:sksol}. Next,
\begin{multline*}
O(\epsilon):\;\;\;\; \left(-a_0^2 \frac{\partial ^2}{\partial x^2} + (r_0^2 + b_0 r_0)\right) \eta_1 = \frac{k_m}{\rho_0}\frac{\partial \eta_0}{\partial
x} \\
+\frac{k_s}{2\rho_0^2}\frac{\partial^2 }{\partial x^2} \eta_0
-(2 r_0 r_1 + b_0 r_1)\eta_0 \defeq  R
\end{multline*}
Substituting $r_0 = s_l^{\pm}(0)$ on the left hand side leads to a
resonance condition for the right hand side term, denoted by $R$. In
particular for a solution $\eta_1$ to exist, $R$ must lie in the
range space of the linear operator
\begin{equation}
\left(-a_0^2 \frac{\partial^2}{\partial x^2} + (r_0^2 + b
r_0)\right).
\end{equation}
For this self-adjoint operator, the range space is the complement of
its null space $\{\sin(\frac{l x}{2})\}$. This gives the resonance
condition as
\begin{equation*}
\<R,\sin(\frac{l x}{2})\> = 0,
\end{equation*}
where $<\cdot,\cdot>$ denotes the standard inner product in
$L^2(0,2\pi)$. This leads to an equation
\begin{multline}\label{eq:reson-temp}
(2r_0+b_0)r_1 =  \frac{l}{4\pi\rho_0}\int_0^{2\pi}k_m(x)\sin(lx) dx \\
 - \frac{l^2}{8\pi\rho_0^2}\int_{0}^{2\pi}k_s(x)\sin^2(\frac{lx}{2}) dx
\end{multline}
For values of $r_0=s_l^{\pm}(0)$, where $s_l^{\pm}(0)$ is given
by~\eqref{eq:sksol}, the equation above leads to an expression for
perturbation in the two eigenvalues.  We denote these perturbations
as $r_1^{\pm}$. For $r_0 = s_l^{+}(0)$, we have from from
Lemma~\ref{lem:s-ell-symmetric} that $b_0>>|2r_0|$ when $l<<l_c$, which happens for every $l$ as $N \to \infty$ (see eq.~\eqref{eq:l-critical}), so that
\begin{align}\label{eq:r1plus}
r_1^+ \approx \frac{l}{4\pi\rho_0 b_0}\int_0^{2\pi}k_m(x)\sin(l x) dx +
O(\frac{1}{N^2}).
\end{align}
Note that we have dropped the second integral on the right hand side
of~\eqref{eq:reson-temp} because $\frac{1}{\rho_0^2} = O(1/N^2)$ for large $N$. For
$r_0 = s_k^{-}(0)$, $2r_0 \approx -2b_0$ for $l << l_c$ and
\begin{align}\label{eq:r1minus}
r_1^- \approx - \frac{l}{4\pi\rho_0 b_0}\int_0^{2\pi}k_m(x)\sin(l x)
dx + O(\frac{1}{N^2}).
\end{align}
Note that
\begin{eqnarray*}
r_1^+ + r_1^- &=& 0.
\end{eqnarray*}
Putting the formulas for the perturbation to the
eigenvalues~\eqref{eq:r1plus} and~\eqref{eq:r1minus}
in~\eqref{eq:perturb}, we get
\begin{align*}
s_l^+(\epsilon) & \approx s_l^+(0) + \epsilon \frac{l}{4\pi b_0
\rho_0}\int_0^{2\pi}k_m(x)\sin(lx) dx + O(\epsilon^2) +
O(\frac{1}{N^2}),\\ s_l^-(\epsilon) & \approx -b_0 - \epsilon
\frac{l}{4\pi b_0 \rho_0}\int_0^{2\pi}k_m(x)\sin(lx) dx +
O(\epsilon^2)+O(\frac{1}{N^2}).
\end{align*}
Since $s_l^+(0) = O(\frac{1}{N^2})$ for $l < l_c$
(Lemma~\ref{lem:s-ell-symmetric}) and $\rho_0=\frac{N}{2\pi}$, the result
follows. \end{proof-theorem}

\end{document}